\documentclass{amsart}
\usepackage[all]{xy}
\usepackage{amsmath,amsfonts,amssymb,amsthm}
\newtheorem{theorem}[equation]{Theorem}
\newtheorem{lemma}[equation]{Lemma}
\newtheorem{proposition}[equation]{Proposition}
\newtheorem{definition}[equation]{Definition}
\newtheorem{corollary}[equation]{Corollary}
\theoremstyle{remark}
\newtheorem{remark}[equation]{Remark}

\numberwithin{equation}{section}
\DeclareMathOperator\Hom{Hom}
\DeclareMathOperator\fn{NF}
\DeclareMathOperator\Id{Id}
\DeclareMathOperator\lcm{lcm}
\DeclareMathOperator\End{End}
\DeclareMathOperator\Obj{Obj}
\newcommand\CP{{\mathcal P}}
\newcommand\CI{{\mathcal I}}
\newcommand\BR{{\mathbb R}}
\newcommand\BZ{{\mathbb Z}}
\newcommand\inv{^{-1}}
\newcommand\bB{{\bf B}}
\newcommand\bI{{\bf I}}
\newcommand\bJ{{\bf J}}
\newcommand\bG{{\bf G}}
\newcommand\bS{{\bf S}}
\newcommand\bW{{\bf W}}
\newcommand\ba{{\bf a}}
\newcommand\bb{{\bf b}}
\newcommand\bc{{\bf c}}
\newcommand\bd{{\bf d}}
\newcommand\bs{{\bf s}}
\newcommand\bt{{\bf t}}
\newcommand\bu{{\bf u}}
\newcommand\bv{{\bf v}}
\newcommand\bw{{\bf w}}
\newcommand\bx{{\bf x}}
\newcommand\by{{\bf y}}
\newcommand\CB{{\mathcal B}}
\newcommand\CC{{\mathcal C}}
\newcommand\CF{{\mathcal F}}
\newcommand\CO{{\mathcal O}}
\newcommand{\lexp}[2]{\kern\scriptspace\vphantom{#2}^{#1}\kern-\scriptspace#2}
\newcommand{\ie}{{\it i.e}}
\newcommand{\eg}{{\it e.g}}
\newcommand{\cf}{{\it cf}}
\newcommand{\cpb}{{\CC(P_\bullet(\Id))}}
\newcommand{\f}[1]{{f\inv(Y_{\ge #1})}}
\newcommand{\dash}{{\hbox{---}}}
\begin{document}
\author{F.~Digne and J.~Michel}
\title[Garside categories]{Garside and locally Garside categories}
\date{December 26th, 2006}
\maketitle
\section{Introduction}
This  text is our  version of (locally)  Garside categories. Our motivation is
the  example  of  section  \ref{B+(I)},  which  we needed in September 2004 to
understand  some Deligne-Lusztig varieties. Since this example works naturally
in  the setting of arbitrary  Coxeter groups, at that  time we only considered
the  general  case  of  categories  which  are  locally  Garside.  Krammer has
independently  introduced the  notion of  (full) Garside categories \cite{K1},
\cite{K2}.

The  things we  added since  2004 are  that we  noticed that it makes sense to
consider  categories which are only left or  right locally Garside, and that a
sufficient  condition to make things work is a Noetherianness property (before
that, we imposed the homogeneity which comes from an additive length). We also
added  a discussion of the relation between  our definitions and the notion of
Garside  categories,  for  which  we  use  the definition introduced by Bessis
\cite{B}.  We define what we call  left Garside categories in this context;
part  of  this  reflects  inspiring  discussions  we  had with Bessis and with
Krammer in april 2006.

The notion of Garside category has recently been used by Bessis \cite{B2}
to obtain deep
theorems  about braid  groups of  complex reflection  groups. We  have a joint
project with David Bessis and Daan Krammer to write a general survey about the
subject.
This text should be taken as our initial contribution to this project.

\section{Locally Garside and Garside categories}

We  adopt  conventions  for  categories  which  are  consistent with those for
monoids  and with those in algebraic topology:  we write $xy$ for the composed
of the morphisms $A\xrightarrow x B$ and $B\xrightarrow y C$. We consider only
small  categories.  The  morphisms  have  a  natural  preorder  given  by left
divisibility:  if $x=yz$, we  say that $y$ is a left divisor
or a left  factor of $x$, which we
denote  $y\preccurlyeq x$. We write $y\prec x$ if in addition $y\ne x$. And we
say  that  $x$  is  a  right  multiple  of  $y$  (we  have  also evidently the
corresponding  notions when exchanging left and right). We will write $x\in C$
to  say that $x$ is a  {\em morphism} in $C$. We  will write $A\in \Obj(C)$ to
say that $A$ is an object of $C$.

\begin{definition}\label{poset noetherien}
We say that a preordered set is {\em Noetherian} if there does not exist any
bounded infinite strictly increasing sequence.
\end{definition}

Notice  that  such  a  set  is  then  a  poset. We say that a category is left
Noetherian if left divisibility induces a Noetherian poset on the morphisms.
This  notion also makes sense (and we will  use it) for a subset of a category.

We  call right lcm (resp.\  left gcd) a least  upper bound (resp.\ largest lower
bound) for left divisibility. An lcm is unique for a poset.

\begin{definition}\label{localement Garside}
A (small) category $C$  is left locally Garside if
\begin{enumerate}
\item\label{Noetherian} It is left Noetherian.
\item\label{simplifiable a gauche}
It has the left cancellation property, \ie.\ $xy=xz$ implies $y=z$ (in other
words, every morphism is an epimorphism).
\item\label{lcm} Two morphisms which have some common right multiple have a 
right lcm.
\end{enumerate}
\end{definition}

A  left locally Garside  monoid is the  monoid of morphisms  of a left locally
Garside category with only one object.

In  a locally  Garside category,  complements are  well defined, that is, when
$x\preccurlyeq  y$ there is a unique $z$,  the {\em complement} to $y$ of $x$,
such that $x=yz$.

We say that a subcategory $C_1$ is stable by complement if $x,y\in C_1$ and
$x\preccurlyeq y$ in $C$ implies that $x\preccurlyeq y$ in $C_1$.

Given a subcategory stable by complement, we say it is stable by lcm if
when $x,y\in C_1$ have an lcm $z$ in $C$ then $z$ is in $C_1$. From the
stability by complement it follows that $z$ is an lcm in $C_1$. The following
lemma is a transposition to categories of a result of \cite{godelle}.
\begin{lemma}\label{subcategory}
Let  $C$ be a  left locally Garside  category; if $C_1$  is a subcategory
stable by complement and by lcm, then $C_1$ is left locally Garside.
\end{lemma}
\begin{proof}
Axioms \ref{Noetherian} and \ref{simplifiable a gauche} for $C_1$ are inherited
from the same  axioms for $C$. Let  us check axiom \ref{lcm}.  If two morphisms
have a common right multiple in $C_1$ they have an lcm in $C$ which is an
lcm in $C_1$ by assumption.
\end{proof}

Before  going  further,  let  us  look  at consequences and equivalent ways of
formulating the third axiom in the setting of posets.

\begin{lemma}\label{lcm  of  family}  In  a  Noetherian  poset  where  any two
elements which have an an upper bound have a least upper bound, the same holds
for any family of elements.
\end{lemma}
\begin{proof}   Let  $\CF$  be  a  bounded  above  subset  of  the  poset.  By
Noetherianity  there exist maximal elements in  $\CF$. As two maximal elements
of  $\CF$ have a  least upper bound  by assumption, they  have to be equal, so
there is a unique maximal element, which is a least upper bound for $\CF$.
\end{proof}

\begin{lemma}\label{gcd  of family} In a Noetherian  poset which has a minimal
element  and where any two elements which have  an an upper bound have a least
upper bound, any family of elements has a largest lower bound.
\end{lemma}
\begin{proof}  Let $\CF$ be the family in the statement. The elements of $\CF$
have  at least one lower  bound, the minimal element  of the poset. Now we can
apply  the previous lemma  to the family  of the lower  bounds of $\CF$, which
have  thus a least upper bound. This least upper bound is smaller than all the
elements of $\CF$, so is a largest lower bound for $\CF$.
\end{proof}

In  the  context  of  left  divisibility,  we  get  that axiom \ref{localement
Garside}(iii) implies
\begin{corollary}\label{axiom iv}
In a left locally Garside category, a family of morphisms with the same source
have a left gcd.
\end{corollary}

There are contexts where the above corollary is equivalent to the third axiom.
\begin{definition}\label{poset artinian}
We say that a poset is {\em Artinian} or if there does not exist
any bounded infinite strictly decreasing sequence.
\end{definition}

We note that for left divisibility any family of morphisms with same source
has always a lower bound, the identity
morphism. It follows from the previous lemmas, using the reverse order, that
\begin{corollary}\label{(iv)=>(iii)}
A category which satisfies axioms \ref{localement Garside}(i) and (ii),
which is left Artinian and such that any two morphisms with same source have
a left gcd, is left locally Garside (\ie.\ satisfies \ref{localement Garside}(iii)).
\end{corollary}

The  notion of right  locally Garside category  is obtained by exchanging left
and right in the definition of a left locally Garside category.

A  locally Garside category is  a category which is  both locally left Garside
and locally right Garside.

\begin{lemma}\label{artinien}
A  subset of a category stable by left divisibility, left Noetherian and which
verifies a weak form of the left cancellation property, that is $xy=x$ implies
$y=1$,  is  also  {\em  right  Artinian}.
\end{lemma}
\begin{proof}
If $a_n$ is a strictly decreasing sequence for the order $\succcurlyeq$ and if
$a_n=b_na_{n+1}$  then $c_n=b_1b_2\ldots  b_n$ is  an increasing  sequence for
$\preccurlyeq$ all terms of which left divide $a_1$. It is strictly increasing
since    $b_1\ldots    b_{n-1}=b_1\ldots    b_n\Rightarrow    b_n=1\Rightarrow
a_n=a_{n+1}$, the first implication by the weak form of left cancellation.
\end{proof}

It  follows  from  \ref{(iv)=>(iii)}  and  \ref{artinien}  that  for a locally
Garside  category we  can replace  axiom \ref{localement  Garside}(iii) by the
existence  of  a  left  gcd  for  all  pairs of morphisms with same source and
similarly on the right.

The three following definitions are adaptations to one-sided Garside
categories of the definitions of \cite{B}. Definition \ref{Garside} is almost
equivalent to \cite[2.5]{B}.
\begin{definition}\label{left Garside}
A left Garside category $C$ is a left locally Garside category such that there
exists an endofunctor $\Phi$ of $C$ and a natural transformation $\Delta$ from
the  identity functor to $\Phi$ such that the set of left divisors of $\Delta$
generate $C$.
\end{definition}

We  denote  by  $A\xrightarrow{\Delta_A}\Phi(A)$  the  natural  transformation
applied to the object $A$; in the above the left divisors of $\Delta$ mean the
divisors of the various $\Delta_A$ as $A$ runs over the objects of $C$.

For right Garside, we change also the direction of the natural
transformation.
\begin{definition}\label{right Garside}
A right Garside category $C$ is a right locally Garside category such that there
exists an endofunctor $\Phi$ of $C$ and a natural transformation $\Delta$ from
$\Phi$ to the  identity functor such that the set of right divisors of $\Delta$
generate $C$.
\end{definition}

Finally, we define Garside:
\begin{definition}\label{Garside}
A  Garside category  $C$ is  a right  and left  Garside category such that the
functor  $\Phi$ for the right  Garside structure is the  inverse of $\Phi$ for
the left Garside structure, and such that the left and right $\Delta$
coincide.
\end{definition}

By  saying that the right and left  $\Delta$ coincide, we mean that $\Delta_A$
for the left Garside structure is the same as $\Delta_{\Phi(A)}$ for the right
Garside structure.

We will show (\cf.\ \ref{garside bilatere}) that a left Garside category which
is  right Noetherian and such that $\Phi$  has an inverse is Garside.

In the case of a Garside monoid identified with
the endomorphisms of a one-object
Garside category, the functor $\Phi $ is the conjugation by the element
$\Delta$ of the monoid.

\section{Germs for locally Garside categories}
We   introduce  a  convenient  technique   for  constructing  locally  Garside
categories  by introducing the notion  of a {\em germ},  which is some kind of
generating  set  for  categories,  and  giving  conditions  on  a germ for the
generated category to be locally Garside. This section is an adaptation in the
context of categories of section 2 of \cite{BDM}; the main technical
difference being that here we assume neither atomicity nor the existence of a
length function: they are replaced by the Noethianness property.

\begin{definition}\label{germe} A germ $(P,\CO)$ is a pair consisting of a set
$\CO$  of objects,  and a  set $P$  of morphisms  (which have  a source  and a
target,  which  are  objects),  with  a  partially defined ``composition'' map
$m:P\times  P\to P$. For $a,b\in P$ we will write `` $ab\in P$''\ to mean that
$m(a,b)$  is defined; and  in this situation  we denote $ab$  for $m(a,b)$; we
abbreviate  $ab\in P$ and $c=ab$ to $c=ab\in  P$. If we denote by $P(A,B)$ the
set  of  morphisms  in  $\CO$  of  source  $A$  and target $B$, we require the
following axioms:
\begin{enumerate}
\item For all $A\in\CO$, there exists $1_A\in P(A,A)$
such that for any $a\in P(B,A)$ (resp.\ any $a'\in P(A,B)$)
$a=a. 1_A\in P$ (resp.\ $a'=1_A.a'\in P$).
\item For $a,b,c \in P$, we have $ab, (ab)c \in P$ if and
only if $bc, a(bc)\in P$ and in this case $a(bc)=(ab)c$.
\end{enumerate}
\end{definition}
We will write 1 instead of $1_A$ when the context makes clear that the source
of this morphism is $A$.
A  {\em path} in $P$  is a sequence of  morphisms $(p_1,\ldots,p_n)$ such that
the  target of $p_i$  is the source  of $p_{i+1}$. If  $(x_1,\ldots,x_n)$ is a
path  such that  for some  bracketing of  this sequence the product $x_1\ldots
x_n$  is defined in  $P$, then by  axiom \ref{germe} (ii)  the product is also
defined,  and has the same value, for  any bracketing of the sequence. We will
denote by $x_1\ldots x_n\in P$ this situation (and $x_1\ldots x_n$ the product
when this situation occurs).
\begin{definition}\label{categorie engendree}
The category generated by the germ $(P,\CO)$ is the category with objects
$\CO$ defined by generators and relations as follows: the generators are $P$,
and the relations are $ab=c$ whenever $c=ab \in P$.
\end{definition}

We  write $C(P,\CO)$, or $C(P)$  when there is no  ambiguity, for the category
generated  by  the  germ  $(P,\CO)$.  We  can  give  an explicit model for the
morphisms  of $C(P)$  in terms  of equivalences  classes of  paths in $P$. The
equivalence   relations  between   paths  is   generated  by   the  elementary
equivalences:
$$(p_1,\ldots,p_{i_1},p_i,p_{i+1},\ldots,p_n)\sim
(p_1,\ldots,p_{i_1},p'p'',p_{i+1},\ldots,p_n)$$
when $p_i=p'p''\in P$, and $(1)\sim ()$.

The  composition of  morphisms in  $C(P)$ is  defined by  the concatenation of
paths. The next lemma shows that this extends the partial product in $P$.
\begin{lemma}
Let $(x_1,\ldots,x_n)$ be a path equivalent to the the single-term path
$(y)$. Then $x_1\ldots x_n\in P$ and $x_1\ldots x_n = y$.
\end{lemma}

\begin{proof}
The assumption implies that there exists a sequence of elementary equivalences
$$l_0=  (x_1,\ldots,x_n)  \sim  l_1  \sim  \cdots  \sim  l_k= (y)$$ where each
equivalence $l_{j-1} \rightarrow l_{j}$ is either
\begin{itemize}
\item a contraction
$$(p_1,\ldots,p_{i-1},p',p'',p_{i+1},\ldots,p_m) \rightarrow
(p_1,\ldots,p_{i-1},p'p'',p_{i+1},\ldots,p_m)$$
\item an expansion
$$(p_1,\ldots,p_{i-1},p'p'',p_{i+1},\ldots,p_m) \rightarrow
(p_1,\ldots,p_{i-1},p',p'',p_{i+1},\ldots,p_m)$$
\end{itemize}
where  $p'p''\in  P$.  We  may  assume  that  $k$ is minimal. If there is any
expansion in the sequence, let $(p_1,\ldots,p_{i-1},p'p'',p_{i+1},\ldots,p_m)
\rightarrow  (p_1,\ldots,p_{i-1},p',p'',p_{i+1},\ldots,p_m)$  be  the  last
one. Since all subsequent steps are contractions we have $m$
subsequent  steps and $y=p_1\ldots  p_{i-1}p'p''p_{i+1}\ldots p_m\in P$. Since
any             bracketing             of             the             sequence
$(p_1,\ldots,p_{i-1},p',p'',p_{i+1},\ldots,p_m)$  has the  same value,  we see
that  we could start  with the bracketing  $\ldots(p'p'')\ldots$, and thus get
from  $(p_1,\ldots,p_{i-1},p'p'',p_{i+1},\ldots,p_m)$ to $(y)$  in $m-1$ steps
instead  of $m+1$ whence a contradiction unless there are only contractions in
a minimal sequence of equivalences leading to $(y)$, whence the result.
\end{proof}
We have the following
\begin{corollary}\label{divise}
$P$ identifies with a subset of $C(P)$ stable by taking left or right factors.
\end{corollary}
\begin{proof}
Indeed, if two morphisms of $P$ are equal in $C(P)$ the above lemma (using the
particular  case $n=1$) shows that they are equal  in $P$. And a left or right
factor  in $C(P)$  of $y$  in  the above  lemma is a product $x_1\ldots x_i$ or
$x_i\ldots x_n$ and is thus in $P$.
\end{proof}

Just  as for a  category, we say  that a germ  $(P,\CO)$ is left Noetherian if
left divisibility induces a Noetherian poset on $P$.

Let  us remark that  a germ with  a superadditive length,  that is, a function
$P\xrightarrow   l\BZ_{\ge   0}$   such   that   $l(ab)\geq   l(a)+l(b)$   and
$l(a)=0\Leftrightarrow a=1$ is left and right Noetherian.

\begin{lemma}\label{1.4}
Let $C$  be a category and $P$ be a set of morphisms which generates $C$.
Let  $X$ be  a set  of morphisms  of $C$  with same  source satisfying
\begin{enumerate}
\item $X$ is stable  by taking left factors,
\item $X$ is a bounded Noetherian poset for left divisibility,
\item If $x\in X$, $y,z\in P$ and $xy,  xz\in X$ then $y$ and $z$ have a common
right  multiple $m$ such that $xm\in X$.
\end{enumerate}
Then  $X$ is the set of left divisors of some morphism of $C$.
\end{lemma}
\begin{proof}
Since  $X$ is a bounded Noetherian poset for $\preccurlyeq$, there exists a maximal element $g\in
X$  for $\preccurlyeq$. Let us  prove by contradiction
that $X$ is  the set of left divisors of
$g$.  First we notice that otherwise
$E=\{x\prec  g\mid \exists  u\in P,  xu\in X,  xu\not\preccurlyeq g\}$  is not
empty: indeed let $y\in X$ be such that $y\not\preccurlyeq g$ and let $x$ be a
maximal  common factor  of $y$  and $g$;  then $x$  in $E$,  since if we write
$y=xu_1\ldots  u_k$  with  $u_i\in  P$  and  $k$  minimal,  then $k\neq 0$ and
$xu_1\neq  x$, thus $xu_1\not\preccurlyeq  g$ (by maximality  of $x$). Let now
$x\in  E$ be maximal for $\preccurlyeq$ and let $u$ be as in the definition of
$E$.  As $x\prec g$, there is $v\in P$ such that $x\prec xv\preccurlyeq g$. As
$xu$  and $xv$ are both in $X$, the assumption on $X$ implies that $u$ and $v$
have a common multiple $m$ such that $xm\in X$. As $xm$ is a right multiple of
$xu$  we have  $xm\not\preccurlyeq g$.  Thus if  $v'$ is  a maximal  such that
$v\preccurlyeq  v'\preccurlyeq  m$  and  $xv'\preccurlyeq  g$, we have $x\prec
xv'\prec  xm$ and $xv'\in E$ (taking for the  $u$ in the definition of $E$ any
element   of  $P$  such  that  $v'u\preccurlyeq  m$),  which  contradicts  the
maximality of $x$.
\end{proof}

\begin{definition}\label{pregarside}
A germ $(P,\CO)$ is left locally Garside if
\begin{enumerate}
\item[(G1)]\label{G1}  It is left Noetherian.
\item[(G2)]\label{G2}  If two morphisms in $P$ have a common right multiple in $P$,
they have a right lcm in $P$.
\item[(G3)]\label{G3} If two morphisms $u,v\in P$ have a right lcm $\Delta_{u,v}\in
P$ and if $x\in P$ is such that $xu, xv\in P$ then $x\Delta_{u,v}\in P$.
\item[(G4)]\label{G4}  For $z\in C(P)$ and  $x,y\in P$, the equality $zx=zy$
implies $x=y$.
\end{enumerate}
\end{definition}

\begin{remark}\label{inject}
We  note that  \ref{G4}(G4) is  the only  axiom which  does not involve only a
check  on elements of $P$. However, in practical applications, it will be easy
to  check since  it is  automatically verified  if there  is an  injective map
compatible  with  multiplication  from  $P$  into  a  category  with  the left
cancellation property.
\end{remark}

We have a weak form of right cancellation
\begin{lemma}\label{simplifiable a droite}
If $P$ is a germ satisfying (G1) and (G4) then the equality $xy=y\in P$ implies $x=1$.
\end{lemma}
\begin{proof}
From $xy=y$ we deduce that for all $n$ we have $x^ny=y$, so $x^n$ is an increasing
sequence for $\preccurlyeq$ which is bounded by $y$ so has to be constant for $n$ large
enough by (G1). But $x^n=x^{n+1}$ implies $x=1$ by (G4).
\end{proof}
We  will show that the category generated by a locally Garside germ is locally
Garside  by directly constructing normal forms  for elements of $C(P)$. We fix
now a locally Garside germ $(P,\CO)$.

\begin{proposition}
\label{pgcd dans P}
Any family of morphisms in $P$ with same source have a left gcd.
\end{proposition}
\begin{proof}
This is a consequence of lemma \ref{gcd of family}, whose assumption is true by \ref{G2}(G2).
\end{proof}
\begin{proposition}
\label{alpha dans P}
If  $x,y \in P $  are such that the  target of $x$ is  the source of $y$,
then there is a unique maximal $z$ such that $z\preccurlyeq y $ and $xz\in P$.
\end{proposition}
\begin{proof}
This  time  we  apply  lemma  \ref{1.4}  to  the  set  $X$  of  $u$  such that
$u\preccurlyeq y$ and $xu\in P$. $X$ inherits left Noetherianity from $P$ thus
it is enough to check that if $u, v,  w\in P$ are such that $uv,uw\in X$ then
they  have a  right lcm  $\Delta_{v,w}$ and  $xu\Delta_{v,w}\in P$ (which will
imply $u\Delta_{v,w}\in X$). As $uv$ and $uw$ are left factors of $y\in P$, by
axiom  \ref{G4}  (G4)  they  have  a  common  multiple,  in $P$ by corollary
\ref{divise},  thus by \ref{G2}  (G2) they have  a lcm $\Delta_{v,w}$ which by
axiom \ref{G3} (G3) satisfies $xu\Delta_{v,w}\in P$.
\end{proof}
\begin{definition}
\label{alpha2 et omega2}
Under   the   assumptions   of   proposition   \ref{alpha   dans   P}  we  set
$\alpha_2(x,y)=xz$  and  we  write  $\omega_2(x,y)$  for the morphism $t\in P$
(unique    by   axiom   \ref{G4}   (G4))   such   that   $y=zt$.   Thus   $xy=
\alpha_2(x,y)\omega_2(x,y)$.
\end{definition}
\begin{proposition}
\label{alpha2(xy,z) et omega2(xy,z)}
For $x,y,z,xy\in P$ we have
\begin{enumerate}
\item $\alpha_2(xy,z)=\alpha_2(x,\alpha_2(y,z))$.
\item $\omega_2(xy,z)=\omega_2(x,\alpha_2(y,z))\omega_2(y,z)$.
\end{enumerate}
\end{proposition}
\begin{proof}
Let   us  show   (i).  Define   $u,v  \in   P$  by   $\alpha_2(xy,z)=xyu$  and
$\alpha_2(y,z)=yv$.  As $yu\in  P$, $u\preccurlyeq  z$, we  have by definition
$yu\preccurlyeq     yv$.     Similarly     $xyu\preccurlyeq    \alpha_2(x,yv)=
\alpha_2(x,\alpha_2(y,z))$;    let    thus    $u'\in    P$    be   such   that
$\alpha_2(x,\alpha_2(y,z))=xyuu'$.  Since  $xyuu'\preccurlyeq  xyv$  by  axiom
\ref{G4}  (G4) have  $uu'\preccurlyeq v$;  as $v\preccurlyeq  z$
we have $uu'\preccurlyeq  z$, and as
$xyuu'\in P$, we  have $u'=1$ by maximality of $u$  in the definition of
$\alpha_2(xy,z)$, which gives (i).

Let us show (ii). Using $xyz=\alpha_2(xy,z)\omega_2(xy,z)$ and
\begin{multline*}
\alpha_2(xy,z)\omega_2(x,\alpha_2(y,z))\omega_2(y,z)=\alpha_2(x,\alpha_2(y,z))
\omega_2(x,\alpha_2(y,z))\omega_2(y,z)=\\
x\alpha_2(y,z)\omega_2(y,z)=xyz,\end{multline*}
which  comes from (i) and definition \ref{alpha2  et omega2}, we will get (ii)
if  we can simplify $\alpha_2(xy,z)$ between  these two expressions for $xyz$.
We  apply axiom \ref{G4} (G4) if we show that both sides of (ii) lie in $P$.
It  is the case for $\omega_2(xy,z)\in P$  by definition, thus we have to show
that   $\omega_2(x\alpha_2(y,z))\omega_2(y,z)\in  P$.   Define  $u\in   P$  by
$\alpha_2(y,z)=yu$,   so   that   $u\omega_2(y,z)=z$,   and   $u_1\in   P$  by
$\alpha_2(x,\alpha_2(y,z))=xyu_1$,                   so                   that
$xyu_1\omega_2(x,\alpha_2(y,z))=x\alpha_2(y,z)=xyu$.                      Then
$u_1\omega_2(x,\alpha_2(y,z))\in   P$   as   it   is   a   right   factor   of
$\alpha_2(y,z)\in   P$,   thus   $u_1\omega_2(x,\alpha_2(y,z))=u$   (by  axiom
\ref{G4}(G4))                                                           thus
$u_1\omega_2(x,\alpha_2(y,z))\omega_2(y,z)=u\omega_2(y,z)=z$   which   implies
that  $\omega_2(x\alpha_2(y,z))\omega_2(y,z)\in P$ as it  is a right factor of
an element of $P$.
\end{proof}

\begin{proposition}\label{alpha}
There  is a unique map  $\alpha:C(P)\to P$ which is  the identity on $P$, such
that  for $x,y\in P$ we have $\alpha(xy)=\alpha_2(x,y)$, and such that for any
$u,v\in C(P)$ we have $\alpha(uv)=\alpha(u\alpha(v))$. In addition $\alpha(u)$
is the unique maximal left factor in $P$ of $u$.
\end{proposition}
\begin{proof}
We  will  define  $\alpha$  on  the  paths  in  $P$,  and  then check that our
definition  is  compatible  with  elementary  equivalence.  As $\alpha$ is the
identity  on $P$,  we need  that $\alpha(())=1$  and that  $\alpha((y))=y$ for
$y\in P$. The conditions we want impose that
\begin{equation}\label{*}
\alpha(p_1,\ldots,p_k)=\alpha_2(p_1,\alpha(p_2,\ldots,p_k)).\end{equation}
By  induction on $k$, this already shows  that $\alpha$ is unique. We will now
show  by  induction  on  $k$  that  $\alpha$ is compatible with the elementary
equivalence    $(p_1,\ldots,p_k)\sim(p_1,\ldots,p_ip_{i+1},\ldots,p_k)$   when
$p_ip_{i+1}\in P$. If this equivalence is applied at a position $i>1$, formula
\ref{*}   shows  that  compatibility   for  paths  of   length  $k-1$  implies
compatibility for paths of length $k$. If $i=1$ we have to compare
$\alpha_2(p_1,\alpha(p_2,\ldots,p_k))$
and $\alpha_2(p_1p_2,\alpha(p_3,\ldots,p_k))$.
But $\alpha_2(p_1p_2,\alpha(p_3,\ldots,p_k))=
\alpha_2(p_1,\alpha_2(p_2,\alpha(p_3,\ldots,p_k)))$,
by \ref{alpha2(xy,z) et omega2(xy,z)} (i) and
$\alpha_2(p_2,\alpha(p_3,\ldots,p_k))=
\alpha(p_2,p_3,\ldots,p_k)$ by \ref{*}, whence the result that
$\alpha$ is well defined by \ref{*} on $C(P)$.

Similarly, if $u=(u_1,\ldots,u_m)$ and $v=(v_1,\ldots,v_n)$, we show
that $\alpha(uv)=\alpha(u\alpha(v))$ by induction on $m+n$. Indeed
\begin{multline*}
\alpha(u_1,\ldots,u_m,v_1,\ldots,v_n)=
\alpha_2(u_1,\alpha(u_2,\ldots,u_m,v_1,\ldots,v_n))=\hfill\\
\hfill\alpha_2(u_1,\alpha(u_2,\ldots,u_m,\alpha(v)))=
\alpha(u_1,u_2,\ldots,u_m,\alpha(v)),
\end{multline*}
by respectively \ref{*}, the induction hypothesis, and \ref{*} again.

Finally  we show that $\alpha(u)$ is the maximal left factor in $P$ of $u$. It
is  by definition an element of $P$ which is  a left factor of $u$. If we have
another expression $u=pv$ with $p\in P$ then
$\alpha(u)=\alpha(pv)=\alpha(p\alpha(v))=\alpha_2(p,\alpha(v))$ so $\alpha(u)$
is a right multiple of $p$.
\end{proof}
\begin{proposition}\label{omega}
There is a unique map $\omega: C(P)\to C(P)$ such that
for $x,y\in P$ we have $\omega(xy)=\omega_2(x,y)$, and such that
for any $u,v\in C(P)$ we have $\omega(uv)=\omega(u\alpha(v))\omega(v)$.
\end{proposition}
\begin{proof}  As in  the previous  proposition we  define $\omega$ on paths by
induction.  We must have $\omega(x)=1$  for $x\in P$ and  for a path of length
$k\ge  2$  we  must have
\begin{equation}\label{**}
\omega(p_1,\ldots,p_k)=\omega_2(p_1,\alpha(p_2\ldots p_k))\omega(p_2,\ldots,p_k)
\end{equation}
This  proves the unicity  of $\omega$, and  again we show  by induction on $k$
that  this is
compatible with elementary equivalence.  Again, we come to the case
of  an elementary  equivalence occurring  in the  first term, \ie., to compare
$\omega(p_1,p_2,\ldots,p_k)$  and $\omega(p_1p_2,p_3,\ldots,p_k)$ when $p_1p_2
\in P$. We have
\begin{multline*}
\omega(p_1,\ldots,p_k)=\omega_2(p_1,\alpha(p_2\ldots p_k))\omega(p_2,\ldots,p_k)
=\hfill\\
\omega_2(p_1,\alpha(p_2\ldots p_k))\omega_2(p_2,\alpha(p_3\ldots
p_k))\omega(p_3,\ldots,p_k)=\hfill\\
\hfill\omega_2(p_1p_2,\alpha(p_3\ldots    p_k))
\omega(p_3,\ldots,p_k)=\omega(p_1p_2,p_3,\ldots,p_k)
\end{multline*}
by  respectively \ref{**},  \ref{**}, \ref{alpha2(xy,z)  et omega2(xy,z)} (ii)
and \ref{**} whence the result.

We show similarly for
$u=(u_l,u_2,\ldots,u_m)$ and $v=(v_1,\ldots,v_n)$ that
$\omega(uv)=\omega(u\alpha(v))\omega(v)$  by induction on $m+n$.
We have
$$\begin{aligned}
\omega(uv)&=\omega(u_1,\ldots,u_m,v_1,\ldots,v_m)\\&=
\omega_2(u_1,\alpha(u_2\ldots
u_mv_1\ldots v_n))\omega(u_2,\ldots,u_m,v_1,\ldots,v_n)\\
&=\omega_2(u_1,\alpha(u_2\ldots u_mv_1\ldots v_n))
\omega(u_2,\ldots,u_m,\alpha(v_1\ldots v_n))\omega(v_1,\ldots,v_n)\\
&=\omega_2(u_1,\alpha(u_2\ldots
u_m\alpha(v_1\ldots
v_n)))\omega(u_2,\ldots,u_m,\alpha(v_1\ldots v_n))\omega(v_1,\ldots,v_n)\\
&=\omega(u_1,u_2,\ldots,u_m,\alpha(v_1\ldots
v_n))\omega(v_1,\ldots,v_n),\end{aligned}$$
by respectively  \ref{**}, the induction hypothesis,
\ref{alpha}  and \ref{**}, whence the result.
\end{proof}

We are now ready to define normal forms for morphisms in $C(P)$.
\begin{definition}\label{normale}
We  call {\em  normal form}  of a  morphism $x\in  C(P), x\ne  1$ a decomposition
$x=x_1\ldots  x_k$ such that $x_i\in P, x_k\ne  1$ and such that for all $i$
we have $x_i=\alpha(x_i\ldots x_k)$.
\end{definition}
We  notice that we have $x_i\ne 1$ for all $i$ since an element $\ne1$ has
a  non-trivial  $\alpha$.  We  declare  that  the  normal form of $1$ is the
trivial decomposition ($k=0$).

We will show the existence of normal forms in \ref{existence de forme normale}
and their unicity in \ref{unicite de la forme normale}.
We first show another characterization.
\begin{proposition}\label{forme normale}
A  decomposition $x_1\ldots x_k$ where $x_i\in P, x_k\ne 1$ is a normal form
if  and only if for all $i\le  k-1$ the decomposition $x_ix_{i+1}$ is a normal
form.
\end{proposition}
\begin{proof}
If   $x_1,\ldots,   x_k$   is   a   normal   form  then  $x_i=\alpha(x_i\ldots
x_k)=\alpha(x_i\alpha(x_{i+1}\ldots       x_k))=\alpha(x_ix_{i+1})$       thus
$x_ix_{i+1}$  is a normal form. Conversely, assuming by decreasing induction
on $i$, that
$\alpha(x_{i+1}\ldots x_k)=x_{i+1}$, we have  $\alpha(x_i\ldots x_k)=
\alpha(x_i\alpha(x_{i+1}\ldots  x_k))=\alpha(x_ix_{i+1})=x_i$.
\end{proof}
The above statement implies that any product of consecutive terms in a normal
forms is itself a normal form.

\begin{proposition}\label{forme normale de xy}
If $x_1\ldots x_k$ is a normal form and $y\in P$, there exist
decompositions $x_i=x'_ix''_i$ such that either
$(yx'_1)(x''_1x'_2)\ldots (x''_{k-1}x'_k)x''_k$, if $x''_k\neq 1$, or
$(yx'_1)(x''_1x'_2)\ldots (x''_{k-1}x'_k)$ otherwise, is a normal form
of $yx_1\ldots x_k$.
\end{proposition}
\begin{proof}
This  is obtained by  recursively applying $\alpha(ab)=\alpha(a\alpha(b))$: we
first     write     $\alpha(yx)=yx'_1$     where     $x_1=x'_1x''_1$,     then
$\alpha(x''_1x_2\ldots x_k)=
\alpha(x''_1\alpha(x_2\ldots x_k))=\alpha(x''_1x_2)=x''_1x'_2$ etc\dots
\end{proof}

\begin{corollary}\label{existence de forme normale}
Normal forms exist.
\end{corollary}
\begin{proof}
We  proceed by induction  on $k$ for  an $x\in C(P)$  of the form $x=p_1\ldots
p_k$  with $p_i\in P$. By induction we may as well assume that $p_2\ldots p_k$
is  a normal  form. The  previous proposition  shows then  how to  construct a
normal form for $x$.
\end{proof}

We now show that $C(P)$ has the let cancellation property. We will deduce
it from the following property of $\omega$.
\begin{proposition}\label{alpha.omega}
If $x\in C(P)$, then $\omega(x)$ is the unique $y\in C(P)$
such that $x=\alpha(x)y$.
\end{proposition}
\begin{proof}  We show that $x=\alpha(x)y$  implies $y=\omega(x)$ by induction
on the number of terms in a normal form of $y$. If $y=1$ then $x\in P$ and the
result holds. The assumptions at step $k$ are now that $x=\alpha(x)y$ and that
some  decomposition  $y=y_1\ldots  y_k$  is  a  normal  form. By the induction
hypothesis and the equality $y=y_1(y_2\ldots y_k)$ we get $\omega(y)=y_2\ldots
y_k$, thus $y=\alpha(y)\omega(y)$.
Thus $\omega(x)=\omega(\alpha(x)y)=\omega(\alpha(x)\alpha(y))\omega(y)$.
On the other hand
$\alpha(x)=\alpha(\alpha(x)y)=\alpha(\alpha(x)\alpha(y))=
\alpha_2(\alpha(x),\alpha(y))$.
Thus $\omega(\alpha(x)\alpha(y))=\omega_2(\alpha(x),\alpha(y))=\alpha(y)$
where the last equality is by definition of
$\omega_2$. Putting things together we get $\omega(x)=\alpha(y)\omega(y)=y$.
\end{proof}
\begin{corollary} \label{simplifiable}
$C(P)$ has the left cancellation property.
\end{corollary}
\begin{proof}
We  want to  show that  for any  $x,y,z\in C(p)$  the equality $xy=xz$ implies
$y=z$.  By induction  on the  number of terms in a decomposition of $x$ into a
product  of  elements  of  $P$,  we  may  assume  that $x\in P$. Define $b$ by
$\alpha(xy)=xb$;  then  $b$  is  unique  since  $P$  has the left cancellation
property.  Let $y'$ be an  element such that $by'=y$,  which is possible since
$b\preccurlyeq\alpha(y)\preccurlyeq  y$.  By  proposition \ref{alpha.omega} we
have $y'=\omega(xy)$. Thus if similarly we define $z'$ as an element such that
$bz'=z$ we have $z'=\omega(xy)=y'$ thus $z=bz'=by'=y$.
\end{proof}
\begin{corollary}\label{unicite de la forme normale}
Normal forms are unique.
\end{corollary}
\begin{proof}
If $x=x_1\ldots x_k$ is a normal form then $x_1=\alpha(x)$ is uniquely defined
and  by proposition  \ref{alpha.omega} we  have $x_2\ldots  x_k=\omega(x)$; we
conclude by induction on $k$.
\end{proof}

For $x\in C(P)$, we denote by $\nu(x)$ the minimum number of terms
in a decomposition of $x$ into a product of elements of $P$.
\begin{lemma}\label{nu}
The normal form of $x$ has $\nu(x)$ terms.
\end{lemma}
\begin{proof}
The  proof is by induction on $\nu(x)$.  We assume the result for $\nu(x)=k-1$
and  will prove  it for  $\nu(x)=k$. Let  then $x=x_1\ldots  x_k$ be a minimal
decomposition  of $x$.  By induction,  the normal  form of $x_2\ldots x_k$ has
$k-1$ terms, so we may as well assume that $x_2\ldots x_k$ is normal. By lemma
\ref{forme  normale de xy} the normal form of $x$ has $k-1$ or $k$ terms. Thus
it has $k$ terms, whence the result.
\end{proof}

\begin{lemma}\label{nu croissant}
If  $x$ is a right factor of $y$ then $\nu(x)\leq \nu(y)$.
\end{lemma}
\begin{proof}
Since an element can be obtained from a right factor by repeatedly multiplying
on the left by elements of $P$, proposition \ref{forme normale de xy} shows that
a right factor has less terms in its normal form.
\end{proof}
If $x\in C(P)$ has normal form $x=x_1\ldots x_n$ we have
$\omega^k(x)=x_{k+1}\ldots x_n$ for $k\geq 0$
(it is 1 if $k\geq n$),  and the $k$-th term of the normal
form of $x$ is $\alpha(\omega^{k-1}(x))$.
\begin{lemma}\label{omega divise}
For $a\in P$ and $x\in C(P)$, if $\omega^k(x)=\omega^k(ax)$ for some $k$ then
$\alpha(\omega^{k-1}(ax))\succcurlyeq\alpha(\omega^{k-1}(x))$.
\end{lemma}
\begin{proof}
In this proof (only) we will still call normal form a product with a certain
number of trailing $1$'s.
Let $x=x_1\ldots x_n$ be a normal form of $x$.
Then by proposition \ref{forme normale de xy}
which is still valid with our present definition of normal forms,
we can write $x_i=x'_ix''_i$ for all $i$, so that
$(ax'_1)(x''_1x'_2)\ldots(x''_{n-1}x'_n)x''_n$ is a normal form of $ax$.
By assumption $(x''_kx'_{k+1})\ldots (x''_{n-1}x'_n)x''_n
=\omega^k(ax)=\omega^k(x)=x_{k+1}\ldots x_n$,
so identifying these two normal forms we get $x''_n=1$ and
$x''_ix'_{i+1}=x_{i+1}=x'_{i+1}x''_{i+1}$ for $i\geq k$, whence by decreasing induction
on $i$ and lemma \ref{simplifiable a droite} we get $x''_i=1$  so that
$x_i=x'_i$ for $i\geq k$. Whence
$\alpha(\omega^{k-1}(ax))=
x''_{k-1}x'_k=x''_{k-1}x_k=x''_{k-1}\alpha(\omega^{k-1}(x))$.
\end{proof}
\begin{proposition}\label{C(P) noetherien}
$C(P)$ is Noetherian for left divisibility.
\end{proposition}
\begin{proof}
We  have to show  that no infinite  sequence $x_1\prec x_2\prec x_3\cdots\prec
x_n\prec\cdots\preccurlyeq   x$   exists;   we   proceed   by   induction   on
$\nu(x)$.
If $\nu(x)=1$ the sequence consists of elements of $P$ which contradicts
the  Noetherianity of  $P$. The  sequence $\alpha(x_i)$  is non-decreasing and
bounded  by $\alpha(x)$ so is  constant at some stage  by the Noetherianity of
$P$.  Truncating the  previous terms  and simplifying  by the  common value
$a_1$ of $\alpha(x_i)$,
we  get  a  an  infinite  sequence  bounded by $a_1\inv x$. If $\nu(a_1\inv
x)<\nu(x)$ then we are done by induction. Otherwise we can repeat the
same argument for another step, introducing the common value $a_2$ of
$\alpha(a_1\inv x_i)$, etc\dots; after $k$ such steps we will still
have $\nu((a_1\ldots a_k)\inv x)=\nu(x)$.
But this implies  by lemma \ref{omega divise}
that $\omega^{\nu(x)-1}((a_1\ldots a_h)\inv x)$
is a decreasing sequence of elements of $P$, so it has to be constant at some
stage. Truncating at
this stage we may assume that the last term of the normal form of
$\nu((a_1\ldots a_h)\inv x)$ is equal to the last term of the normal
form of $x$.
Lemma \ref{omega divise}  gives then that $\alpha(\omega^{\nu(x)-2}((a_1\ldots a_h)\inv
x))$ is decreasing for right divisibility so has to be constant at
some stage. Truncating again we can assume that in the whole process
the last two terms of the normal form of $(a_1\ldots a_h)\inv x$
are constant. Going on we come to a point where  $(a_1\ldots a_h)\inv x$
itself is constant which means, again by \ref{omega divise},
that $a_h=1$ for $h$ large enough.
\end{proof}

\begin{proposition}\label{Delta_st dans C(P)}
If two elements of $P$ have a common right multiple in $C(P)$ then they have
a right lcm in $P$ (which is also their lcm in $C(P)$.
\end{proposition}
\begin{proof}
We  first observe that if $u,v$ have a common multiple $x\in C(P)$ they have a
common  multiple in $P$, which is  $\alpha(x)$. We may then apply \ref{G2}(G2)
to conclude.
\end{proof}

\begin{proposition}\label{ppcm dans C(P)}
Any family of elements of $C(P)$  who have a common right multiple has
a right lcm. If the family is a subset of $P$ then the lcm is in $P$.
\end{proposition}
\begin{proof}
Assume $(x_i)_{i\in I}$ have a common right multiple. We apply lemma \ref{1.4}
to the set $X$ of elements which divide all the common multiples of the $x_i$.
It  inherits Noetherianity from $C(P)$, and if $u,v \in P$ and $x\in C(P)$ are
such that $xu,xv\in X$ then any right multiple of $xu$ and $xv$ is of the form
$xz$  where $u$ and $v$  divide $z$ (by the  cancellation property in $C(P)$);
thus   the  right  lcm  $\Delta_{u,v}$  of   $u$  and  $v$  (which  exists  by
\ref{Delta_st  dans C(P)}) divides $z$. Thus $x\Delta_{u,v}\in X$. We may thus
apply  \ref{1.4} and the elements  of $X$ are the  divisors of an element which
must be the lcm of the $x_i$. The second statement comes from the fact that if
the $x_i$ are in $P$ and divide $x$ then they divide $\alpha(x)$
and so does their lcm.
\end{proof}

\begin{proposition}\label{pgcd dans C(P)}
Any family of morphisms in  $C(P)$ has a left gcd.
\end{proposition}
\begin{proof}
This is a consequence of lemma \ref{gcd of family}, whose assumption is true by \ref{ppcm dans C(P)}.
\end{proof}

At this stage we have proved the following
\begin{theorem}\label{categorie<-germe}
If the germ $P$ is left locally Garside, so is $C(P)$.
\end{theorem}

This has a converse:
\begin{theorem}\label{categorie->germe}
Let $C$ be a left locally Garside category $C$. Then $(P,\CO)$, where $\CO$ is
the  set of objects of  $C$ and where $P$  is a set of  morphisms of $C$ which
generate  $C$, stable by taking left factors  and right factors, and stable by
taking  right lcm when  they exist, is  a left locally  Garside germ; for this
germ, $C=C(P)$.
\end{theorem}
\begin{proof}
We  first check that  a $P$ such  as above is  a germ. Axiom  (i) of a germ is
clear.  Axiom (ii) (``associativity'') holds for a  set of morphisms as soon as
they  are stable by  taking left and  right factors. The  axioms for a locally
Garside germ \ref{pregarside} are immediate except perhaps axiom \ref{G3} (G3)
(for  \ref{G4}(G4) see  remark \ref{inject}).  If $u,v\in  P$ have a right lcm
$\Delta_{u,v}$ and if $x\in P$ is such that $xu,xv \in P$ then $x\Delta_{u,v}$
is the right lcm in $C$ of $xu$ and $xv$ thus is in $P$.

All   the  relations  of  $C(P)$   hold  in  $C$,  thus   we  have  a  functor
$C(P)\xrightarrow i C$ which is clearly surjective since $P$ generates $C$. We
have to see that $i$ is injective.

Let  us define a function  $\alpha:C\to P$ defined for  $x\in C$ by taking the
largest  (for left divisibility) factor of $x$  in $P$; this exists, since $P$
is stable by right lcm. The formula $\alpha(xy)=\alpha(x\alpha(y))$ holds when
$x\in  P$ since  $\alpha(xy)$ is  by definition  of the  form $xa\in  P$ where
$xa\preccurlyeq  xy$ thus $a\preccurlyeq y$  by the cancellation property thus
$a\preccurlyeq \alpha(y)$.

To    see   that   $i$   is   injective    it   is   enough   to   show   that
$i\circ\alpha=\alpha\circ i$; indeed, by induction on $\nu(x)$ for $x\in C(P)$
we     have     $i(x)=i(y)\Rightarrow     \alpha(i(x))=\alpha(i(y))\Rightarrow
i(\alpha(x))=i(\alpha(y))\Rightarrow  \alpha(x)=\alpha(y)$, the  last equality
since  $i$ is injective on $P$. By  the left cancellation property in $C$ this
implies $i(\omega(x))=i(\omega(y))$ and we conclude by induction.

Let  us show  that for  any $x\in  C(P)$ we  have $i\circ\alpha(x)=\alpha\circ
i(x)$ by induction on $\nu(x)$. Let $x_1\ldots x_n$ be the normal form of $x$.
Then       $\alpha(i(x_1\ldots      x_n))=\alpha(i(x_1)i(x_2\ldots      x_n))=
\alpha(i(x_1)\alpha(i(x_2\ldots  x_n)))=\alpha(i(x_1)i(x_2))$  where  the last
equality  is by  the induction  hypothesis. We  are thus  reduced to  the case
$n=2$,   \ie.\  to  show  that  if   $x,y\in  P$  and  $\alpha_2(x,y)=x$,  then
$\alpha(xy)=x$ in $C$. But this is clear by the definitions of $\alpha$ in $C$
and $\alpha_2$ in $P$.
\end{proof}

\subsection*{Subgerms, and fixed points}
\begin{definition}
If  $(P,\CO)$ is a  germ, we call  {\em subgerm} of  $P$ a pair $(P_1,\CO_1)$
obtained by taking a part $\CO_1$ of the objects $\CO$ and a part $P_1$ of the
morphisms   between  objects  in  $\CO_1$  which  is  stable  by  the  partial
multiplication in $P$, and contains the morphisms $1_A$ for $A\in\CO_1$.
\end{definition}
It is  straightforward to check  that a subgerm is  a germ. Note,  however, that
left  divisibility might  be  quite  different in  $P_1$:  it  is possible  that
$a,ab\in P_1$ but $b\in  P-P_1$ in which case we do  not have $a\preccurlyeq ab$
in $P_1$. If the divisibility in $P_1$ is the restriction of the divisibility in
$P$, we say that $P_1$ is stable by  complement. As in the case of categories we
say that a  subgerm $P_1$ stable by complement  is stable by lcm if  for any two
morphisms in $P_1$ which have a common multiple in $P_1$, their lcm in $P$ is in
$P_1$ (so is an lcm in $P_1$).
\begin{lemma}\label{sous-germe}
If $P$ is a left locally Garside  germ, and $P_1$ a subgerm stable by complement
and lcm, then $P_1$ is left locally Garside.
\end{lemma}
\begin{proof}
Axiom  \ref{G1}(G1) is clearly inherited from $P$ to $P_1$. Axiom \ref{G4}(G4)
is  also inherited from $P$, using  the natural functor $C(P_1)\to C(P)$ which
is  injective  on  $P_1$  (since  its  restriction  to  $P_1$ restricts to the
injection $P_1\to P$, because $P\to C(P)$ is injective).

Let  us check \ref{G2}(G2). If  two elements have a  common multiple in $P_1$,
they  have an lcm in $P$ since $P$ is locally Garside. That lcm is in $P_1$
and is an lcm in $P_1$ by assumption.

Let  us check \ref{G3}(G3). If two elements $u,v$ have a lcm $\Delta_{u,v}$ in
$P_1$, then by assumption their lcm in $P$ is in $P_1$, and must thus
be  equal to $\Delta_{u,v}$. Thus if  $xu,xv\in P_1$ then $x\Delta_{u,v}\in P$
thus $x\Delta_{u,v}\in P_1$ since $P_1$ is stable by partial multiplication.
\end{proof}

\begin{lemma}\label{stable alpha2}
If,  under the assumptions of \ref{sous-germe}, in addition $P_1$ is stable by
$\alpha_2$ (that  is  when  $x,y\in  P_1$  then  $\alpha_2(x,y)\in P_1$), then
$C(P_1)$ injects in $C(P)$.
\end{lemma}
\begin{proof}
We have to show that the natural functor $C(P_1)\xrightarrow i C(P)$ which sends
a  path to  the  corresponding path  is  injective. Since  $i$  is injective  on
$P_1$  it is  enough to  show  that $i$  preserves  normal forms;  by the  local
characterization of normal forms it is enough to show that the image of a 2-term
normal form is a  normal form. But that is a consequence of  the fact that $P_1$
is stable by $\alpha_2$: indeed, if $(x,y)$ is a 2-term normal form in $P_1$ and
$\alpha_2(x,y)=xz$ in $P$ then by assumption we have $xz\in  P_1$, whence $z\in
P_1$ and $z$ divides $y$ in $P_1$, as  $P_1$ is stable by complement, so that $z=1$
since $(x,y)$ is normal in $P_1$.
\end{proof}
\begin{proposition}\label{points fixes}
Let  $P$ be a left locally Garside germ and let $\sigma$ be an automorphism
of  $C(P)$  stabilizing  $P$;  let  $P^\sigma$  (resp.\  $C(P)^\sigma$)  be the
subgerm  (resp.\ the subcategory) of the $\sigma$-fixed morphisms and objects
of  $P$ (resp.\ $C(P)$);  then $P^\sigma$ is  a left locally  Garside germ and
$C(P)^\sigma=C(P^\sigma)$.
\end{proposition}
\begin{proof}
The  unicity of complement, lcm and  $\alpha_2$ shows that  $P^\sigma$ is stable by 
complement, lcm
and  $\alpha_2$.  Thus $P^\sigma$  is a  left locally  Garside germ and the
natural  functor  $C(P^\sigma)\xrightarrow  i  C(P)^\sigma$  is injective. As,
given  a $\sigma$-fixed morphism of $C(P)$, all
the terms of its normal form are in $P^\sigma$ by the unicity of normal forms,
we get that $i$ is surjective.
\end{proof}

\subsection*{A counterexample}
We give an example to show that the endomorphisms of an object in a locally
Garside category are not necessarily a locally Garside monoid.

Let  $(P,\CO)$ be the germ where $\CO=\{X,Y\}$ and where there are seven
morphisms:
$s,t\in\End(Y)$,     two    elements     $a,b\in\Hom(X,Y)$,    two    elements
$u,v\in\Hom(Y,X)$, plus one additional morphism resulting from the only composition
defined in $P$, given by $as=bt$. The axioms of a germ as well as \ref{G2}(G2)
and  \ref{G3}(G3) are easy (associativity and \ref{G3} (G3) are empty, and the
only morphisms having a common right multiple are $a$ and $b$, and this multiple is
unique;  the same holds on  the left for $s$  and $t$). Let us prove \ref{G4}
(G4) and its right analogue. There is an additive length on $C(P)$ defined by
$l(a)=l(b)=l(s)=l(t)=1$.  If $zx=zy$ for $x,y\in P$  and $z\in C(P)$, then $x$
and  $y$ have the same length. If this length is $2$ then $x=y$ since there is
only  one element of length $2$ in $P$.  Assume then the length is $1$. If $x$
is  neither $s$ nor $t$,  no relation in $C(P)$  for the word $zx$ can involve
$x$,  thus $x=y$.  It remains  to consider  the case  $x=s$ and $y=t$, \ie.\ an
equality  $zs=zt$. If $z$ has no decomposition ending by $a$ again no relation
can involve the terminal $s$. If $z=z'a$, no relation in $C(P)$ can change the
terminal  $a$ into another morphism,  in particular into $b$.  Thus the word $z'at$
cannot  be changed into $z'as$ so we are finished. A similar reasoning applies
on the right, which finishes the proof that $C(P)$ is locally Garside.

On the other hand, in $\End(X)$ the morphisms $a$ and $b$ have two minimal-length
common right multiples $asu=btu$ and $asv=btv$, thus no lcm.

\section{Atoms}
We  call {\em  atom} a  morphism in  a category  (resp.\ a germ) which does not
admit any proper right or left factor.

Note  that if the  category is left  Noetherian and has  the left cancellation
property,  then by \ref{simplifiable a droite} having no proper left factor is
equivalent to having no proper right factor.

We  say that a germ (resp.\ a category)  is {\em atomic} if any morphism in the
germ (resp.\ category) is a product of atoms.

In an atomic category $C$, a set $P$ of morphisms generates $C$ if and only if
it contains the atoms of $C$.

\begin{proposition} \label{engendre}
A  category $C$ which has the left cancellation property and is right and left
Noetherian (\eg.\ a locally Garside category) is atomic.
\end{proposition}
\begin{proof}
Let  us show that any morphism $f\in C$ is a product of atoms. By the analogue
of   \ref{artinien}   on   the   right,   which   is   applicable   thanks  to
\ref{simplifiable a droite}, we know that $C$ is left Artinian, thus $f$ has a
left  factor which has no proper left factor, and which is thus an atom by the
remark above.

Let  thus $a_1$ be an atom which is a left factor of $f$ and let $f_1$ be such
that  $f=a_1f_1$. We may similarly write $f_1=a_2  f_2$ where $a_2$ is an atom
which  is a left factor of $f_1$, etc\dots and if $f$ is not equal to a finite
product  $a_1\ldots  a_n$  we  would  get  an  infinite increasing sequence of
factors of $f$ which would contradict left Noetherianity.
\end{proof}
\begin{remark}
The  cancellation  property  is  necessary  in  the  the above proposition. An
example  of a left and  right Noetherian monoid without  atoms is given by the
set $\{x_i\}_{i\in\BZ_{\ge0}}\cup\{x_\infty\}$
where $x_i$ is the infinite  sequence
beginning by $i$ times $1$ followed by
all  $0$s  and $x_\infty$  is the infinite sequence with
all terms equal to 1; the product is by term-wise multiplication.
We  have $x_ix_j=x_{\inf(i,j)}$.
\end{remark}
\begin{proposition}
Under  the assumptions of \ref{points fixes}, if in addition $C(P)$ is atomic,
then  $C(P^\sigma)$ is also atomic with atoms the right lcm of orbits of atoms
of  $P$ (for  the orbits  which have  a common  multiple) which  are not right
multiples of another such lcm.
\end{proposition}
\begin{proof}  Let us first see that any element $x\in C(P^\sigma)$ is
divisible on the left by such an lcm.
Let  $s\in P$ be an atom such that $s \preccurlyeq x$. Then any element of the
orbit  of $s$ also divides  $x$, thus their lcm $\overline s$,  which is in
$P^\sigma$, also does divide $x$. By the left cancellation property, if
we write $x=\overline s x_1$, then we also have $x_1\in C(P)^\sigma$; we can
apply the same process to $x_1$ to get an $x_2$, etc\dots. By \ref{artinien}
the sequence $x_i$ is finite thus $x$ is a finite product of such lcm.
Finally, such an lcm which
is not divisible by another is clearly an atom in $C(P)^\sigma$.
\end{proof}

We will now give conditions in terms of atoms which imply the properties for
locally Garside.

\begin{proposition} A Noetherian atomic germ $P$ satisfying \ref{G4}
(G4) satisfies (G2) and  (G3) if and only if
\begin{enumerate}
\item[(G2$'$)]\label{G2'}  If two atoms have a right common multiple
in $P$ then they have a right lcm in $P$.
\item[(G3$'$)]\label{G3'}  If two atoms  $s$ and  $t$ have a right lcm
$\Delta_{s,t}\in P$  and $xs,xt\in P$  then
$x\Delta_{s,t}\in P$.
\end{enumerate}
\end{proposition}
\begin{proof}
These conditions are necessary, as a special case of \ref{G2}(G2) and
\ref{G3}(G3).

Let  us show that (G2$'$) implies (G2).  Assume $x,y\in P$ have a common right
multiple  in $P$. We apply  lemma \ref{1.4} to the  set $X$ of elements of $P$
which  are left factors of  all common right multiples  of $x$ and $y$, taking
for  the $P$ of  \ref{1.4} the atoms  in $P$. We  may do so since $X$ inherits
Noetherianity  from $P$, and the assumption  of \ref{1.4} comes from (G4) and
(G2$'$):  if $z\in  X$ and  $s$ and  $t$ are  atoms such that $zs,zt\in X$, by
(G4)  $s$ and $t$ have  a common right multiple,  thus by (G2$'$) they have a
right  lcm  $\Delta_{s,t}$  and  $z\Delta_{s,t}=\lcm(zs,zt)$  is  in $X$. The
common multiple of elements of $X$ given by \ref{1.4} is the desired lcm.

We  now show (G3). Let $u,v,x\in P$ be such  that $u$ and $v$ have a right lcm
$\Delta_{u,v}\in  P$ and such that $xu,xv\in  P$. This time we apply \ref{1.4}
to  $X=\{y\in P\mid xy\in P  \text{ and } y\preccurlyeq\Delta_{u,v}\}$, taking
again  for the  $P$ of  \ref{1.4} the  atoms. Again $X$ inherits Noetherianity
from  $P$; we have $u,v\in X$ by assumption.  Assume now that $y\in X$ and the
atoms  $s$  and  $t$  are  such  that  $ys,yt  \in X$, \ie.\ $xys,xyt\in P$ and
$ys\preccurlyeq\Delta_{u,v}, yt\preccurlyeq\Delta_{u,v}$. By (G4) $s$ and $t$
have   a  common  right  multiple  $P$,  thus  by  (G2$'$)  they  have  a  lcm
$\Delta_{s,t}\in    P$   and   by    (G3$'$)   $xy\Delta_{s,t}\in   P$.   Thus
$y\Delta_{s,t}\in  X$  and  the  assumption  of  \ref{1.4}  holds.  The common
multiple   of  the  elements   of  $X$  given   by  \ref{1.4}  is  necessarily
$\Delta_{u,v}$   since  it   is  a   multiple  of   both  $u$  and  $v$.  Thus
$\Delta_{u,v}\in X$ which implies $x\Delta_{u,v}\in P$.
\end{proof}
\section{Garside categories}

In  a left Garside category $C$ the  smallest set containing the left divisors
of  $\Delta$ and stable by taking left  and right factors forms a left locally
Garside germ $P$ such that $C=C(P)$ by \ref{categorie->germe}. The elements of
$P$ are called the {\em simples} of the category.

\begin{remark}
For  a left locally Garside  category $C$, we could  call {\em simples of $C$}
the  set of morphisms of a chosen germ.  Note that if $C$ is in addition right
Noetherian,  there exists always a minimal such  set, which is the minimal set
$P$  of morphisms  of $C$  stable by  taking left  and right factors and right
lcm's  and generating  $C$; indeed  this set  exists and  is unique, since $C$
itself  has these properties and an intersection of sets with these properties
is  also a set  with these properties  by \ref{engendre} (the only non-trivial
property  to check  for an  intersection is  that it  still generates $C$; but
\ref{engendre}  shows that a subset  generates $C$ if and  only if it contains
the atoms, which is a condition stable by intersection). Then $P$ is a locally
Garside germ and $C=C(P)$ by \ref{categorie->germe}. \end{remark}

The set of simples in a left Garside category is not necessarily minimal in the
sense above.

A  simple  $f$  has  a  {\em  complement  to $\Delta$} denoted $\tilde f$, and
defined  by  $f\tilde  f=\Delta$  (it  is  unique  by  the  left  cancellation
property).  If  $f\in  C$,  as  $\Delta$  is a natural transformation from the
identity   to  $\Phi$  we  have  $f\Delta=\Delta\Phi(f)$  whence,  using  left
cancellation  by  $f$,  we  get  $\Delta=\tilde  f\Phi(f)$, which can also be
written $\tilde{\tilde f}=\Phi(f)$.

If  $\Phi$ is an  automorphism this shows  that the set  of left factors of
$\Delta$  is the same as the set of  right factors of $\Delta$. 

\begin{remark}
In  a left  Garside category,  a set  $P$ of  morphisms stable  by taking left
factors and complements to $\Delta$ is stable by taking right factors. Indeed,
if   $ab$  and   $a$  are   in  $P$   then  $ab(ab)\tilde{}=a\tilde   a$  thus
$b(ab)\tilde{}=\tilde a$ thus $b\in P$ as a left factor of $\tilde a$.
\end{remark}

\begin{proposition}\label{x divise Delta n}
In a left Garside category, left divisibility makes the
set of morphisms with same source into a lattice.
\end{proposition}
\begin{proof}
It  is enough to show that any two morphisms with the same source have a right
lcm. If they are simple, they divide $\Delta$ so we are done. Otherwise, given
$x\in  C$, we show by induction  on $n$ that $x\preccurlyeq\Delta^n$ where $n$
is  the number of terms of the normal form of $x$ with respect to $P$. Indeed,
if   $x=x_1\ldots   x_n$   is   the   normal  form,  by  induction  $x_2\ldots
x_n\preccurlyeq      \Delta^{n-1}$     thus     $x_1\ldots     x_n\preccurlyeq
x_1\Delta^{n-1}$.     But     by     definition     of    $\Phi$    we    have
$x_1\Delta^{n-1}=\Delta^{n-1}\Phi^{n-1}(x_1) \preccurlyeq \Delta^n$.
\end{proof}

\begin{proposition}\label{garside bilatere}
A left Garside category which is right Noetherian and such that $\Phi$ is an
automorphism is Garside.
\end{proposition}
\begin{proof}
We first show that such a category  $C$ has the right cancellation property.
Indeed, if $xa=ya$ we have seen in the proof of \ref{x divise Delta n}
that $a\preccurlyeq\Delta^n$ for some $n$
whence $x\Delta^n=y\Delta^n \Leftrightarrow \Delta^n\Phi^n(x)=\Delta^n\Phi^n(y)$
which implies by left cancellation that $x=y$.

We then observe that for two simples $f,g$ we have $g\preccurlyeq f
\Leftrightarrow \tilde g \succcurlyeq \tilde f$ (the implication from left to
right uses the left cancellation property and from right to left the right
cancellation property). This implies that a left gcd of $f$ and $g$ transports
by $\tilde{}$ to a left lcm of $\tilde f$ and $\tilde g$, and conversely a
right lcm transports to a right gcd.

We can argue similarly for arbitrary morphisms by considering the complement
to a suitable $\Delta^n$ instead of the complement to $\Delta$.

We thus get that $C$ is right locally Garside. Since as remarked above, the
fact that $\Phi$ is an autoequivalence implies that the right divisors of
$\Delta$ are the same as the left divisors, the category is right Garside
for the same $\Delta$ and $\Phi\inv$, so we are done.
\end{proof}

The following proposition points to a possible alternative definition of left
Garside categories.
\begin{proposition}\label{lcm=>garside}
A   left   locally   Garside   category   which   has   a   germ   $P$  as  in
\ref{categorie->germe} such that the morphisms in $P$ with a given source have
a right lcm is left Garside.
\end{proposition}
\begin{proof}
Given an object $A$, we define $\Delta_A$ to be the right lcm of the morphisms
of source $A$. The elements of $P$ are the left divisors of the $\Delta _A$ so
these  divisors generate the  category. On the  morphisms of $P$  we define an
operation  $f\mapsto\tilde  f$  by  the  equality  $f\tilde  f=\Delta$,  using
cancellation.  We then define a functor $\Phi$ which maps $A$ to the target of
$\Delta_A$  and  a  map  f  to  $\tilde{\tilde  f}$.  To show that $\Phi$ is a
functor,  since any morphism is a composition of elements of $P$, it is enough
to check that it is compatible with partial composition.

For  $f,g,fg\in P$ we have $fg(fg)\tilde{\phantom a}=\Delta=f\tilde f$ so that
$g(fg)\tilde{\phantom     a}=\tilde    f$,     whence    $g(fg)\tilde{\phantom
a}\tilde{\tilde   f}=\Delta=g\tilde  g$.  By  left  cancellation,  this  gives
$(fg)\tilde{\phantom a}\tilde{\tilde f}=\tilde g$, whence $(fg)\tilde{\phantom
a}\tilde{\tilde          f}\tilde{\tilde         g}=\Delta=(fg)\tilde{\phantom
a}(fg)\tilde{\tilde{\phantom   a}}$   and   by   cancellation   $\tilde{\tilde
f}\tilde{\tilde g}=(fg)\tilde{\tilde{\phantom a}}$ which is what we wanted.

Finally we note that the equality $f\Delta=\Delta\Phi(f)$ for $f\in P$ extends
to  the same  equality for  arbitrary maps  in the  category, which shows that
$\Delta$  is  indeed  a  natural  transformation  from the identity functor to
$\Phi$.
\end{proof}
\section{The conjugacy category}

Conjugation in a monoid or a category is defined as: $w$ is conjugate to $w'$
if there exists $x$ such that $xw'=wx$. In a category, this condition implies
that $w$ and $w'$ are endomorphisms of some object.
\begin{definition}
Given  a category  $C$, the  {\em conjugacy  category} of  $C$ is the category
whose objects are the endomorphisms of $C$ and where $\Hom(w,w')= \{x\in C\mid
wx=xw'\}$.
\end{definition}
We can extend this definition to simultaneous conjugation of a family of
elements, to get the simultaneous conjugacy category.
If $wx=xw'$ as in the above definition we will write $w'=w^x$ and $w=\lexp xw'$.
\begin{proposition}
If  $C$ is a left (resp.\  right) locally Garside category, its (simultaneous)
conjugacy category is also. Further, one can take as simples for the conjugacy
category the morphisms which are induced by simples  of $C$.
\end{proposition}
\begin{proof}
Let us denote by $\CC$ the conjugacy category of $C$.
Since $\CC$ clearly inherits Noetherianity and cancellability
from  $C$, we have just to show the existence of lcm for morphisms which have a
common  multiple. We will  actually show that  lcm and gcd in $C$ of morphisms
of the conjugacy category are lcm and gcd in the conjugacy category.

We  can rephrase  the condition  $x\in\Hom_\CC(w,?)$ as
$x\preccurlyeq  wx$. If we look at simultaneous conjugation of a family $\CF$,
it  will be the simultaneous condition
$x\preccurlyeq wx$  for all $w\in\CF$. Suppose  that $x,y,w\in C$ are  such that
$x\preccurlyeq wx$ and that $xy\preccurlyeq  wxy$; define $w'$ by $xw'=wx$; then
by cancellation $y\preccurlyeq  w'y$, so that $x$ and $xy$  in $\CC$ imply $y\in
\CC$. Suppose  now that $x\preccurlyeq wx$  and $y\preccurlyeq wy$, and  that $x$
and $y$ have a  right lcm $z$ in $C$. Then using  the left cancellation property
we see that $wz$  is the right lcm of $wx$ and $wy$  thus $x\preccurlyeq wz$ and
$y\preccurlyeq  wz$  from  which  it  follows  that  $z\preccurlyeq  wz$,  \ie.\
$z\in\Hom_\CC(w,?)$, so is  the right lcm of  $x$ and $y$ in $\CC$  by the first
part of the proof.

Similarly  the condition  $x\in\Hom_\CC(?,w)$  can  be written  $xw\succcurlyeq
x$,  and   if  $x,y\in\Hom_\CC(?,w)$  have   a  left   lcm  $z$  we   get  that
$z\in\Hom_\CC(?,w)$ and is the left lcm of $x$ and $y$ in $\CC$.

The   second   assertion   of   the  proposition,   follows   from   the   fact 
that    if    $x\preccurlyeq   wx$    then    $\alpha(x)\preccurlyeq\alpha(wx)= 
\alpha(w\alpha(x))\preccurlyeq  w\alpha(w)$  which   shows  that  $\alpha(x)\in 
\Hom_\CC(w,?)$ (and similarly on the right).
\end{proof}
The following is a straightforward consequence of the proposition:
\begin{corollary}\label{germ for CC}
If $P$ is a germ for $C$ and
if we take the germ for the conjugacy category as in the above proposition,
then the normal form of a morphism in the conjugacy category of $C$ is identical
to its normal form in $C$.
\end{corollary}

\subsection*{\label{B+(I)}The locally Garside category $B^+(\CI)$}

The  locally  Garside  category  that  we  will consider in this subsection is
related to the study of the normalizer of the submonoid generated by a part of
the atoms in an Artin monoid, which has been done by Paris and Godelle.

Let  $(W,S)$ be  a Coxeter  system, and  let $(B^+,\bS)$  be the corresponding
Artin  monoid. Recall that $B^+$
is a locally Garside monoid, with germ the canonical lift
$\bW$ of $W$ in $B^+$ consisting of the elements whose length with
respect to $\bS$ is equal to the length of their image in $W$ with respect to
$S$ (see \eg., \cite{michel}).
Let $\bI_0\subset\bS$ and  let $\CI$ be the  set of conjugates of
$\bI_0$.  Since conjugacy preserves the length (measured with respect to the
generating set $\bS$), we see that any element of $\CI$ is also a subset of
$\bS$. Let $\CC_\CI$ be the connected component of the
(simultaneous) conjugacy category whose objects are $\CI$.
As the monoid $B^+$ is locally
Garside the category $\CC(\CI)$ is left locally Garside.

We denote by $B^+_\bI$ the submonoid of $B^+$ generated by a set
$\bI\subset\bS$.

We recall some definitions and results from \cite{DMR}.

\begin{proposition}
\begin{enumerate}
\item Any $\bb\in B^+$ has a maximal left divisor in $B^+_\bI$, denoted
$\alpha_\bI(\bb)$. We denote by $\omega_\bI(\bb)$ the unique element such
that $\bb=\alpha_\bI(\bb)\omega_\bI(\bb)$.
\item
Let $\bI,\bJ\subset\CI$ and $\bb\in B^+$; 
then $\lexp\bb{B^+_\bJ}\subset B^+_\bI$
if and only if $\lexp{\omega_{\bI}(\bb)}\bJ=\bI$.
\item\label{alpha(b1b2} Let $\bI,\bJ\subset\CI$ and let $\bb_1, \bb_2\in B^+$
be such that $\lexp{\bb_1}\bJ=\bI$ and $\alpha_\bI(\bb_1)=1$; then
$\alpha_\bJ(\bb_2)=1$ if and only if $\alpha_\bI(\bb_1\bb_2)=1$.
\item\label{alphalcm}  If  $\bb_1,  \bb_2\in  B^+$ satisfy
$\bI^{\bb_1}\subset\bS$, $\bI^{\bb_2}\subset\bS$ and $\alpha_\bI(\bb_1)=
\alpha_\bI(\bb_2)=1$,    then   their   right   lcm   $\bc$   also   satisfies
$\alpha_I(\bc)=1$
\end{enumerate}
\end{proposition}
\begin{proof}
(i) is \cite[2.1.5]{DMR}(ii) and (ii) results from \cite[2.3.10]{DMR}.
Let  us prove  (iii). For  $\bs\in \bI$  there exists  $\bs'\in \bJ$ such that
$\bs\bb_1=\bb_1\bs'$.  This element  is then  a common  multiple of  $\bs$ and
$\bb_1$ and has to be their lcm since $\bs'$ is an atom of
$B^+$. So $\bs\preccurlyeq \bb_1\bb_2$ if and only if 
$\bb_1\bs'\preccurlyeq\bb_1\bb_2$, \ie,
$\bs'\preccurlyeq\bb_2$ whence the equivalence of $\alpha_\bJ(\bb_2)=1$
and $\alpha_\bI(\bb_1\bb_2)=1$.

To prove (iv) we will actually show the stronger statement that if in $B^+$ we
have $\bb\preccurlyeq\bc$, $\bI^\bb\subset\bS$  and $\alpha_\bI(\bb)=1$
then $\bb\preccurlyeq\omega_\bI(\bc)$. We proceed by induction on the length
of $\alpha_\bI(\bc)$. If $\alpha_\bI(\bc)=1$ the result is trivial. Otherwise
there exists $s\in\bI$, $s\preccurlyeq\alpha_\bI(\bc)$. By the assumption
$\bI^\bb\subset\bS$ there exists $\bt\in\bS$ such that $\bs\bb=\bb\bt$ or
equivalently $\bs\inv\bb=\bb\bt\inv$. If we write $\bc=\bb\ba$ with $\ba\in
B^+$ then by assumption $\bs\inv\bc=\bb\bt\inv\ba$ is positive \ie, we have
$\bb\bt\inv=\bx\ba\inv$ for some $\bx\in B^+$.
As $\bb\bt\inv=\bs\inv\bb\not\in B^+$, we have 
$\bb\not\succcurlyeq\bt$, whence by unicity of irreducible fractions (see \cite[3.2]{michel})
$\bx=\bb\by$ and $\ba=\bt\by$ for some $\by\in B^+$;
thus $\bt\inv\ba\in B^+$,  \ie. $\bb \preccurlyeq\bs\inv\bc$. We then conclude 
by induction on the length of $\alpha_\bI(\bc)$ that
$\bb\preccurlyeq\omega_\bI(\bs\inv\bc)=\omega_\bI(\bc)$.
\end{proof}

Statement (ii) in the above proposition is a motivation  for restricting the 
next definition to
elements such that $\alpha_\bI(\bb)=1$ (we ``lose nothing'' by doing so).

The following definition makes sense by \ref{alpha(b1b2}
\begin{definition}
We define $B^+(\CI)$ as the category whose set of objects is $\CI$ and such
that the morphisms from $\bI$ to $\bJ$ are the elements $\bb\in B^+$ such that
$\bI^\bb=\bJ$  and $\alpha_\bI(\bb)=1$  (such a  morphism will  be denoted
$(\bI,\bb,\bJ)$).
\end{definition}

By \ref{alpha(b1b2} and \ref{alphalcm} the subcategory $B^+(\CI)$ of $\CC(\CI)$
satisfies the assumptions of lemma \ref{subcategory} and similarly on the
right, so it is locally Garside.

%

We now get a germ for $B^+(\CI)$ from the germ $\bW$ of
the locally Garside monoid $B^+$.
By  \ref{germ for  CC} we have  a germ  $P$ for 
$\CC(\CI)$ consisting of the elements of  $\bW$ which are in $\CC(\CI)$.

\begin{proposition}  Let $\bb$ be a morphism of ${B^+(\CI)}$;
then all the terms of the
normal form in $\CC(\CI)$ of $\bb$ are in ${B^+(\CI)}$.
\end{proposition}
\begin{proof} Let $\bb=\bw_1\ldots\bw_k$ be the normal form of 
$\bb\in\Hom_{B^+(\CI)}(\bI,\bJ)$
in $\CC(\CI)$ (\ie.\ in $B^+$). As $\bw_i\in\CC(\CI)$, we have
$\bI_i=\lexp{\bw_{i+1}\ldots\bw_k}\bJ\subset \bS$  for  all $i$.
Now, as
$\lexp{\bw_1\ldots\bw_{i-1}}\alpha_{\bI_i}(\bw_i\ldots\bw_k)\subset B^+_\bJ$,
so divides $\alpha_\bI(\bb)$, this element has to be 1, whence the result.
\end{proof}
\begin{corollary}
The set of $(\bI,\bw,\bJ)$ in $\CC(\CI)$ such that $\bw\in\bW$ and $\alpha_\bI(\bw)=1$ is a germ
for $B^+(\CI)$.
\end{corollary}

We now identify the germ of the above corollary with a germ constructed in $W$.
It will be convenient to work with  roots instead of subsets of the generators. 
We use the  standard geometric realization of  $W$ as a reflection  group in an 
$\BR$-vector  space $V$  endowed  with a  basis $\Pi$  in  bijection with  $S$. 
The set of roots,  denoted  by $\Phi$  is  the  set  $W\Pi$.  We denote  by  $\Phi^+$ 
(resp.\  $\Phi^-$)  the  elements  of  $\Phi$  which  are  linear  combinations 
with  positive (resp.\  negative) coefficients  of $\Pi$;  a basic 
property  is  that  $\Phi=\Phi^+\coprod\Phi^-$. 
For  $\alpha\in \Pi$  let $s_\alpha$  be the  corresponding element  of $S$ (a
reflection with root $\alpha$).
For  $I\subset\Pi$  we  denote  by  $W_I$ the  subgroup  of  $W$  generated  by 
$\{s_\alpha\mid \alpha\in I\}$; we say that $I$ is spherical if $W_I$ is finite 
and  we then  denote  by  $w_I$ its  longest  element.  A subset  $I\subset\Pi$ 
corresponds to a subset $\bI\subset\bS$. We denote by the same letter the conjugacy class $\CI$
and the corresponding orbit of subsets of $\Pi$.
We say that $w\in W$ is $I$-reduced if 
$w\inv  I\in\Phi^+$.  Being $I$-reduced  corresponds  to  the lift  $\bw\in\bW$ 
having $\alpha_\bI(\bw)=1$. So the germ $P$ identifies with
the set of $(I,w,J)$ such that $I,J\in\CI$ and $wJ=I$. The product
$(I,w,J)(J,w',K)$ is defined in $P$ if and only if $l(ww')=l(w)+l(w')$
and is then equal to $(I,ww',K)$.

We now describe the atoms of $B^+(\CI)$ using the results of \cite{BH}.
If   $I$  is   a  subset   of   $\Pi$  and   $\alpha\in  \Pi$   is  such   that 
$\Phi_{I\cup\{\alpha\}}-\Phi_I$  is finite,  then by  \cite{BH} there  exists a 
unique $v(\alpha,I)\in W_{I\cup\{\alpha\}}$  such that $v(\alpha,I)(\Phi_{I\cup 
\{\alpha\}}^+-\Phi_I^+)\subseteq\Phi_{I\cup\{\alpha\}}^-$ and 
$J=v(\alpha,I)I\subseteq  I\cup\{\alpha\}$; when  $I\cup\{\alpha\}$ is  spherical 
then $v(\alpha,I)= w_{I\cup\{\alpha\}}w_I$.
We have $(J,v(\alpha,I),I)\in P$ when $I\in\CI$
\begin{proposition}
The atoms of $P$ are the elements
$(J,v(\alpha,I),I)$ for $I\in\CI$ and $\alpha\in\Pi-I$ such that
$\Phi_{\{\alpha\}\cup I}-\Phi_I$ is finite.
\end{proposition}
\begin{proof}
By \cite[3.2]{BH} the elements $(J,v(\alpha,I),I)$
as in the proposition generate the monoid.
They are atoms because by \cite[4.1]{BH} the lcm of two such elements,
when it exists, has length strictly larger than either of them.
\end{proof}

\subsection*{The spherical case}

We show now  that $B^+(\CI)$ is Garside when $\bW$  is finite. We recall
that in that case $B^+$ is a Garside monoid, with $\bw_\bS$ as $\Delta$.
We  denote  by $\bs\mapsto\bar\bs$  the  involution  on $\bS$  given  by
$\bs\mapsto\lexp{\bw_\bS}\bs$. This extends  naturally to involutions on
$\CI$ and on $B^+$ that we denote in the same way. We define the functor
$\Phi$  by $\Phi(\bI)=\bar\bI$  and
$\Phi((\bJ,\bw,\bI))=(\bar\bJ,\bar\bw,\bar\bI)$.
The   natural  transformation   $\Delta$  is   given  by   the collection of
morphisms $(\bJ,\bw_\bJ\inv\bw_\bS,\bar\bJ)$.
The  properties
which must be satisfied by $\Delta$ and $\Phi$ are easily checked.

\section{A result \`a la Deligne for locally Garside categories}

In this section we prove a simply connectedness property for the
decompositions into simples for a map in a locally Garside category. This
result is similar (but weaker, see the remark after \ref{Deligne})
to Deligne's result in \cite{Deligne}, but the proof is much simpler and the
result is sufficient for the applications that we have in mind.
The present proof follows a suggestion by Serge Bouc to use a version of
\cite[lemma 6]{Bouc}.

Let $P$ a left locally Garside germ and fix $g\in C(P)$ with $g\neq 1$.
We denote by $E(g)$ the set of decompositions of $g$ into a product of
elements of $P$ different from 1.

Then $E(g)$  is a poset, the order being defined by
$$(g_1,\ldots,g_{i-1},g_i,g_{i+1},\ldots,g_n)>
(g_1,\ldots,g_{i-1},a,b,g_{i+1},\ldots,g_n)$$ if $ab=g_i \in P$.

We recall the definition of the notion of  homotopy in a poset
(which is nothing but a translation of the notion of homotopy
in a simplicial complex isomorphic to $E$ as a poset).
A path from $x_1$  to $x_k$  in $E$  is a sequence
$x_1\ldots x_k$ where each  $x_i$ is  comparable to $x_{i+1}$. The composition
of paths is defined by concatenation. We denote homotopy by $\sim$. It 
is the finest equivalence relation on paths compatible with concatenation and
generated by the two following elementary relations:
$xyz\sim xz$ if $x\le y\le  z$ and $xyx\sim
x$ (resp.\ $yxy\sim y$) when $x\le  y$. Homotopy classes form a groupoid,
as the composition of a paths with source $x$ and of the inverse path is the
constant path at $x$. For $x\in E$ we denote by $\Pi_1(E,x)$ the fundamental
group of $E$ with base point $x$, which is the group of homotopy classes of
loops starting from $x$.

A poset $E$ is said to be {\em simply connected} if it is connected
(there is a path linking any two elements of $E$) and if the fundamental group
with some (or any) base point is trivial.

Note that a poset with a smallest or largest element $x$ is simply connected
since any path $(x,y,z,t,\ldots,x)$ is homotopic to $(x,y,x,z,x,t,x,\ldots,x)$
which is homotopic to the trivial loop.

\begin{theorem}
\label{Deligne} (Deligne)
The set $E(g)$ is simply connected.
\end{theorem}
In fact Deligne, in his more specific setting,
proves the stronger result that $E(g)$ is contractible.
\begin{proof}
First we prove a version of a lemma from \cite{Bouc} on order preserving maps
between posets. For a poset $E$ we put
$E_{\ge  x}=\{x'\in  E\mid  x'\ge x\}$,  which is a simply connected subposet
of $E$ since it has a smallest element. If 
$f:X\to   Y$ is an order preserving map it is compatible with homotopy
(it corresponds to a continuous map between simplicial complexes),
so it induces a homomorphism
$f^*:\Pi_1(X,x)\to \Pi_1(Y,f(x))$.

\begin{lemma} \label{bouc} (Bouc) Let $f:X\to Y$ an order preserving map
between two posets. We assume that $Y$ is connected and that for any $y\in Y$
the poset $\f y$ is connected and non empty. Then $f^*$  is surjective.  If
moreover $\f  y$ is simply connected for all $y$ then 
$f^*$ is  an isomorphism.
\end{lemma}
\begin{proof}
Let us first show that $X$ is connected. Let
$x,x'\in X$; we choose a path $y_0\ldots y_n$  in  $Y$ from
$y_0=f(x)$ to  $y_n=f(x')$.
For $i=0,\ldots,n$, we choose  $x_i\in\f{y_i}$  with $x_0=x$ and  $x_n=x'$.
Then if $y_i\geq y_{i+1}$ we have $\f{y_i}\subset \f{y_{i+1}}$ so that there
exists a path in $\f{y_{i+1}}$  from  $x_i$  to  $x_{i+1}$;  otherwise
$y_i<y_{i+1}$, which implies  $\f{y_i}\supset \f{y_{i+1}}$ and there exists a
path in $\f{y_i}$ from $x_i$ to $x_{i+1}$. Concatenating these paths gives
a path connecting $x$ and $x'$.

We fix now
$x_0\in   X$. Let $y_0=f(x_0)$. We prove that
$f^*:\Pi_1(X,x_0)\to\Pi_1(Y,y_0)$ is surjective. Let
$y_0y_1\ldots  y_n$ with  $y_n=y_0$  be a loop in $Y$.  We lift arbitrarily
this loop into a loop $x_0\dash\cdots\dash x_n$  in
$X$  as above,  (where $x_i\dash x_{i+1}$ stands for a path
from  $x_i$  to  $x_{i+1}$ which is either
in $\f{y_i}$ or in $\f{y_{i+1}}$. Then the path  $f(x_0\dash
x_1\dash\cdots\dash x_n)$  is homotopic to $y_0\ldots y_n$; this can be seen
by induction: let us assume that $f(x_0\dash x_1\cdots\dash x_i)$ is homotopic
to $y_0\ldots y_if(x_i)$;  then the same property holds for $i+1$: indeed
$y_iy_{i+1}\sim  y_if(x_i)y_{i+1}$  as they are two paths in a simply
connected set which is either $Y_{\ge  y_i}$ or
$Y_{\ge y_{i+1}}$; similarly we have $f(x_i)y_{i+1}f(x_{i+1}) \sim f(x_i\dash
x_{i+1})$. Putting things together gives
$$
\begin{aligned}
y_0\ldots y_iy_{i+1}f(x_{i+1})&\sim y_0y_1\ldots y_if(x_i)y_{i+1}f(x_{i+1})\\
&\sim f(x_0\dash\cdots \dash x_i)y_{i+1}f(x_{i+1})\\
&\sim f(x_0\dash\cdots \dash x_i\dash x_{i+1}).
\end{aligned}
$$

We now  prove injectivity  of $f^*$  when all $\f{y}$  are simply  connected.

We first  prove that  if $x_0\dash \cdots\dash  x_n$ and  $x'_0\dash \cdots\dash
x'_n$  are two  loops  lifting the  same  loop $y_0\ldots  y_n$,  then they  are
homotopic.  Indeed,  we get  by  induction  on  $i$ that  $x_0\dash  \cdots\dash
x_i\dash x'_i$  and $x'_0\dash\cdots\dash x'_i$  are homotopic paths,  using the
fact that  $x_{i-1}$, $x_i$, $x'_{i-1}$  and $x'_i$ are  all in the  same simply
connected sub-poset, namely either $\f{y_{i-1}}$ or $\f{y_i}$.

It remains to prove  that we can lift homotopies, which amounts  to show that if
if  we lift  as above  two loops  which differ  by an  elementary homotopy,  the
liftings are  homotopic. If $yy'y\sim y$  is an elementary homotopy  with $y<y'$
(resp.\ $y>y'$), then $\f{y'}\subset\f{y}$ (resp.\ $\f{y}\subset\f{y'}$) and the
lifting of  $yy'y$ constructed as  above is in  $\f{y}$ (resp.\ $\f{y'}$)  so is
homotopic to the trivial path. If  $y<y'<y''$, a lifting of $yy'y''$ constructed
as above  is in $\f{y}$  so is homotopic  to any path  in $\f{y}$ with  the same
endpoints.
\end{proof}

We now prove \ref{Deligne}. By \ref{artinien}  $C(P)$ is right Artinian. Thus if
\ref{Deligne} is  not true there exists  $g\in C(P)$ which is  minimal for right
divisibility such  that $E(g)$ is  not simply connected. Let  $T$ be the  set of
elements  of $P$  which  are left  divisors  of $g$.  By  \ref{ppcm dans  C(P)},
for  any $I\subset  T$  the elements  of  $I$  have an  lcm  $\Delta_I$. We  put
$E_I(g)=\{(g_1,\ldots,g_n)\in E(g)\mid  \forall s\in I,\,  s\preccurlyeq g_1\}$.
The set $E_I(g)$  is the set of  decompositions of $g$ whose first  term is left
divisible by $\Delta_I$.

We  claim  that  $E_I(g)$  is  simply connected  for  $I\neq\emptyset$.  In  the
following,  if   $a\preccurlyeq  b$,  we   denote  by  $a\inv  b$   the  element
$c$  such   that  $b=ac$.  We  apply   \ref{bouc}  to  the  map   $f:  E_I(g)\to
E(\Delta_I\inv   g)$   defined    by   $$(g_1,\ldots,g_n)\mapsto   \begin{cases}
(g_2,\ldots,g_n)&\text{if    $g_1=\Delta_I$}\\   (\Delta\inv_I    g_1,g_2,\ldots
g_n)&\text{otherwise}\\ \end{cases}.$$ This map preserves  the order and any set
$\f{(g_1,\ldots,g_n)}$ has a  least element, namely $(\Delta_I,g_1,\ldots,g_n)$,
so is  simply connected. As by  minimality of $g$ $E(\Delta\inv_I  g)$ is simply
connected \ref{bouc} implies that $E_I(g)$ is simply connected as claimed.

We  now  apply  \ref{bouc}  to  the map  $f:E(g)\to  Y=\CP(T)  -  \{\emptyset\}$
defined  by $(g_1,\ldots,g_n)\mapsto  \{s\in T\mid  s\preccurlyeq g_1\}$,  where
$\CP(T)$  is  ordered   by  inclusion.  This  map  is   order  preserving  since
$(g_1,\ldots,g_n)<(g'_1,\ldots,g'_n)$  implies $g_1\preccurlyeq  g'_1$. We  have
$\f{I}=E_I(g)$, so this set is  simply connected Since $\CP(T) - \{\emptyset\}$,
having a greatest  element, is simply connected \ref{bouc} gives  that $E(g)$ is
simply connected, whence the theorem.
\end{proof}

\section{The categories associated to $P^n$}
Let $P$ be a left locally Garside germ. For any positive integer $n$ we define
a  germ $P_n$ whose objects are the paths of length $n$ in $P$ and such that a
morphism    $\ba\xrightarrow   f   \bb$   where   $\ba=(a_1,\ldots,a_n)$   and
$\bb=(b_1,\ldots,b_n)$,  is given  by a  sequence $f_i$  for $1\le  i\le n+1$,
where  $f_i$ is a morphism from the source  of $a_i$ to the source of $b_i$ for
$i\le n$ and $f_{n+1}$ is a morphism from the target of $a_n$ to the target of
$b_n$, with the additional condition that $f_i\preccurlyeq a_i$ for $i\le
n$  and that, if we define $f'_i$ by $a_i=f_if'_i$, using left cancellability,
we then have $b_i=f'_if_{i+1}$ for $1\le i\le n$.

The     composition    of    two    morphisms    $(a_1,\ldots,a_n)\xrightarrow
f(b_1,\ldots,b_n)\xrightarrow  g  (c_1,\ldots,c_n)$  is  defined in $P_n$ when
$f_i  g_i\preccurlyeq a_i$. We  then set $(fg)_i=f_ig_i$,  which satisfies the
conditions   for  being   in  $P_n$;   indeed,  if   we  define  $(fg)'_i$  by
$a_i=(fg)_i(fg)'_i$   the   equality   to   prove  $c_i=(fg)'_i(fg)_{i+1}$  is
equivalent  by left cancellation to  $(fg)_i c_i=a_i(fg)_{i+1}$, which is true
since        $f_ig_ic_i=f_ig_ig'_ig_{i+1}=f_ib_ig_{i+1}=f_if'_if_{i+1}g_{i+1}=
a_if_{i+1}g_{i+1}$.

Divisibility in $P_n$ is then given by the following result:
\begin{lemma}\label{divisible}
The  morphism  $\ba\xrightarrow  f  \bb$  left  divides  in $P_n$ the morphism
$\ba\xrightarrow  g\bc$ if and only if  $f_i\preccurlyeq g_i$. Then there is a
unique morphism $\bb\xrightarrow h\bc$ such that $fh=g$, where $h_i$ is given
by $f_ih_i=g_i$, using the left cancellation property.
\end{lemma}
\begin{proof}
By the description of the product it is clear that if $h$ is such that
$fh=g$ then $f_ih_i=g_i$. Let us see that conversely this implies that
$h$ is a morphism from $\bb$ to $\bc$. Indeed
$g_ic_i=a_ig_{i+1}\Leftrightarrow f_ih_ic_i=f_if'_if_{i+1}h_{i+1}$ which
implies $h_ic_i=b_ih_{i+1}$, so if we define $h'_i$ by $b_i=h_ih'_i$ we get
$c_i=h'_ih_{i+1}$ as wanted.
\end{proof}
\begin{lemma} $P_n$ is a germ. \end{lemma}
\begin{proof}
Axiom  \ref{germe}  (i)  is  clear,  the  identity morphism being given by the
sequence $f_i=1$.

Let  us  check  axiom  \ref{germe}  (ii).  Let  us  consider  three  morphisms
$\ba\xrightarrow  f\bb\xrightarrow  g\bc\xrightarrow  h\bd$.  From the
definition of the product in $P_n$, 
since when they are defined, we have
$(fg)_i=f_ig_i$ and $(gh)_i=g_ih_i$, the condition for
$(fg)h\in P_n$ and for $f(gh)\in P_n$ is the same, namely that
$f_ig_ih_i\preccurlyeq
a_i$, and both products are defined by the sequence $f_ig_ih_i$ so are equal.
\end{proof}

We will also consider the two subgerms of $P_n$ defined by one of the two
additional conditions:
\begin{definition}
\begin{enumerate}
\item\label{a'=1} The subgerm $P_n(\Id)$ has the same objects as $P_n$ and its
morphisms verify the additional condition $f_1=f_{n+1}=1$.
\item\label{F}
Let  $F$  be  a  functor  from  $C(P)$  to  itself. The objects of the subgerm
$P_n(F)$ are the paths $(a_1,\ldots,a_n)$ such that the target of $a_n$ is the
image  by $F$ of  the source of  $a_1$, and the  morphisms of $P_n(F)$ are the
morphisms of $P_n$ verifying the condition $f_{n+1}=F(f_1)$.
\end{enumerate}
\end{definition}
Note  that the condition on the objects of  $P_n(F)$ is such that they have an
identity  morphism. A connected  component of the  category $C(P_n(\Id))$ is a
``category  of  decompositions''  of  a  given  morphism  in  $C(P)$,  while a
connected component of $C(P_n(F))$ corresponds to a connected component of the
``category  of $F$-twisted  conjugacy'' for  $C(P)$. Note that $P_n(\Id)$ is a
subgerm  of $P_n$ stable by  taking left and right  factors, while $P_n(F)$ is
not.

The germ $P_n(\Id)$ was inspired by a conversation with Daan Krammer.
The germ $P_n(F)$ mimics the ``divided categories'' of David Bessis.

Since
$(fg)_i=f_ig_i$  we  can  extend  the  map  $f\mapsto  f_i:P_n\to  P$ to a map
$f\mapsto  f_i: C(P_n)\to C(P)$; this corresponds to the ``product of the
$i$-th column'', equal to $f_i=f_{1i}f_{2i}\ldots f_{ni}$,
in the following picture of a map in $C(P_n)$.
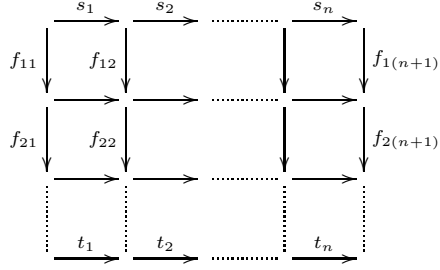
$$\xymatrix{
\ar[r]^{s_1}\ar[d]_{f_{11}}
&\ar[r]^{s_2}\ar[d]_{f_{12}}&\ar@{.}[r]&\ar[r]^{s_n}\ar[d]&\ar[d]^{f_{1(n+1)}}\\
\ar[r]\ar[d]_{f_{21}}&\ar[r]\ar[d]_{f_{22}}&\ar@{.}[r]&\ar[r]\ar[d]
&\ar[d]^{f_{2(n+1)}}\\
\ar[r]\ar@{.}[d]&\ar[r]\ar@{.}[d]&\ar@{.}[r]&\ar[r]\ar@{.}[d]&\ar@{.}[d]\\
\ar[r]^{t_1}&\ar[r]^{t_2}&\ar@{.}[r]&\ar[r]^{t_n}&\\
}
$$
\begin{theorem}\label{C_n(P) localement Garside}
$P_n$ is left locally Garside.
\end{theorem}
\begin{proof}
Let  us check  Noetherianity (\ref{G1}  (G1)). Let  us consider  an increasing
sequence  $(g_k)$ of morphisms all dividing a  morphism $f$ from $(a_1,\ldots,a_n)$ to
$(b_1,\ldots,b_n)$. By lemma 
\ref{divisible} this increasing sequence corresponds to
an  increasing  sequence  of  left  factors  of  $f_i$  for  each $i$ . By the
Noetherianity  of $P$ each of these sequence becomes  constant at some stage
so $g_k$ itself becomes constant and we are done.

We  now  check  left  cancellability  (\ref{G4}  (G4)). Assume that we have an
equality    $fg=fh$   where    $f\in   C(P_n)$    and   $g,h\in   P_n$.   Then
$f_ig_i=(fg)_i=(fh)_i=f_ih_i$  for  all  $i$,  and  by  left cancellability in
$C(P)$ we deduce $g_i=h_i$ for all $i$ q.e.d.

We  now check axiom \ref{G2} (G2). If $f$ and $g$ have a common right multiple
$h$  in $P_n$, then by lemma \ref{divisible} for all $i$ the morphism $h_i$ is
a right multiple of $f_i$ and $g_i$, so $f_i$ and $g_i$ have a right lcm $k_i$
in   $P$.  From  $f_i\preccurlyeq  a_i$   and  $g_i\preccurlyeq  a_i$  we  get
$k_i\preccurlyeq  a_i$, so  $k_i$ defines  a morphism  $k$ in  $P_n$ which is
clearly an lcm for $f$ and $g$.

The  axiom \ref{G3} (G3)  can be similarly  deduced from the corresponding
axiom in $P$.
\end{proof}

\begin{theorem}
\begin{enumerate}
\item\label{Cn(P,F)}
If $F$ preserves right lcms, the category $C(P_n(F))$ is left Garside.
\item\label{Cn(P)}
If $C(P)$ is left Garside, then  $C(P_n)$ also.
\item\label{Cn(P,Id)}
$C(P_n(\Id))$ is  left Garside.
\end{enumerate}
\end{theorem}
\begin{proof}  
We  first check that  the above categories  are left locally  Garside. We have
seen  this for $C(P_n)$ in  \ref{C_n(P) localement Garside}. For $C(P_n(\Id))$
and  $C(P_n(F))$,  since  $P_n(\Id)$  and  $P_n(F)$  are  subgerms of the left
locally  Garside germ $P_n$,  by lemma \ref{sous-germe}  we have just to check
that  they are stable by right complement and lcm. The stability by complement
is  clear  from  the  formula  $(fg)_i=f_ig_i$  for  a  product, since then if
$(fg)_i$  and $f_i$ satisfy the condition for $P_n(\Id)$ (resp. $P_n(F)$) then
$g_i$  will  also  satisfy  it.  Similarly,  since by the proof of \ref{C_n(P)
localement  Garside} the lcm of  $f$ and $g$ is  obtained by taking the lcm of
$f_i$  and $g_i$, it will obviously  satisfy the condition for $P_n(\Id)$, and
also for $P_n(F)$ using the assumption that $F$ preserves lcms.

Thus,  by \ref{lcm=>garside},  we just  have to  check that  in each  of these
categories  the morphisms in the germ with a given source have a right lcm.

For  $C(P_n)$, let $\Delta_A$ be the  natural transformation starting from the
object  $A$ corresponding  to the  left Garside  structure on  $C(P)$, and let
$\ba=(a_1,\ldots,a_n)$  be an object of $P_n$. Then $f_i=a_i$ for $i\le n$ and
$f_{n+1}=\Delta_A$,  where $A$ is the target  of $a_n$, defines a morphism $f$
from  $\ba$ in $P_n$. This morphism is  clearly multiple of any other morphism
from $\ba$.

For the category $C(P_n(F))$, we take the morphism given by $f_i=a_i$ for
$i\le n$ and $f_{n+1}=F(a_1)$. It is clear that it is a multiple in $P_n$
of any morphism from $(a_1,\ldots,a_n)$ which is in in $P_n(F)$; it is also
clear that the quotient is in $P_n(F)$ since $g_{n+1}h_{n+1}=F(g_1h_1)$
and $g_{n+1}=F(g_1)$ imply $h_{n+1}=F(h_1)$ by cancellation.

Finally for $C(P_n(\Id))$, we define by induction for $i\geq 2$ morphisms
$f'_{i-1}$    and   $f_i$   by   the   rules   $a_{i-1}=f_{i-1}f'_{i-1}$   and
$f'_{i-1}f_i=\alpha(f'_{i-1}a_i)$. If we have another morphism $g$ from $\ba$,
we see by the same induction that $g_i\preccurlyeq f_i$ and $g'_i\succcurlyeq
f'_i$.
\end{proof}
Let us spell out the value of $\Phi$ in the first two categories.

For  $C(P_n)$,  we  have  $\Phi((a_1,\ldots,a_n))=  (a_2,\ldots,a_n,\Delta_A)$
where $A$ is the target of $a_n$. If $\ba\xrightarrow f\bb$ is given by $f_i$,
we  have  $\Phi(f)_i=f_{i+1}$  for  $i\le n$ and $\Phi(f)_{n+1}=\Psi(f_{n+1})$
where  $\Psi$ is  the endofunctor  corresponding to  the assumed  left Garside
structure on $C(P)$.

In         $C(P_n(F))$,         if         $\ba=(a_1,\ldots,a_n)$         then
$\Phi(\ba)=(a_2,\ldots,a_n,F(a_1))$;  and $\ba\xrightarrow  f\bb$ is  given by
$f_i$, we have $\Phi(f)_i=f_{i+1}$ for $i\le n$ and $\Phi(f)_{n+1}=F(f_2)$.

\subsection*{The case of a right locally Garside $P$}
In the case where $P$ is right and left locally Garside, we can compute the 
normal form of a morphism and will deduce that a morphism $f$ is determined
by the $f_i$.
\begin{lemma} 
Let $f$ and $g$ be morphisms in $P_n$ such that the target
of $f$ is the source of $g$; then $\alpha(fg)$ is the
morphism whose $i$-th component $\alpha(fg)_i$ is the left gcd of
$f_ig_i$ and $s_i$ for $i=1,\ldots,n$ and $\alpha(fg)_{n+1}=
\alpha(f_{n+1}g_{n+1})$.
\end{lemma}
\begin{proof}
We have the following commutative diagram
$$
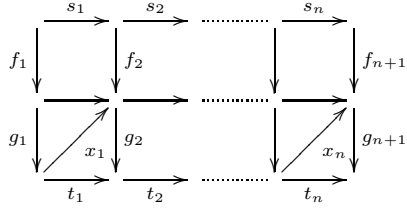
\xymatrix{
\ar[r]^{s_1}\ar[d]_{f_1}
&\ar[r]^{s_2}\ar[d]^{f_2}&\ar@{.}[r]&\ar[r]^{s_n}\ar[d]&\ar[d]^{f_{n+1}}\\
\ar[r]\ar[d]_{g_1}&\ar[r]\ar[d]^{g_2}&\ar@{.}[r]&\ar[r]\ar[d]
&\ar[d]^{g_{n+1}}\\
\ar[r]_{t_1}\ar[ur]_{x_1}&\ar[r]_{t_2}&\ar@{.}[r]&\ar[r]_{t_n}\ar[ur]_{x_n}&\\
}
$$
For $i=1,\ldots,n$, let $f_ia_i$ be the left gcd of $f_ig_i$ and $s_i$ and let
$f_{n+1}a_{n+1}=\alpha(f_{n+1}g_{n+1})$;
we put $g_i=a_ib_i$ for $i=1\ldots,n+1$
and $s_i=f_ia_iv_i$ for $i=1,\ldots,n$. Let $w_i$ be the morphism
$x_ia_i$. The diagram is the following
$$
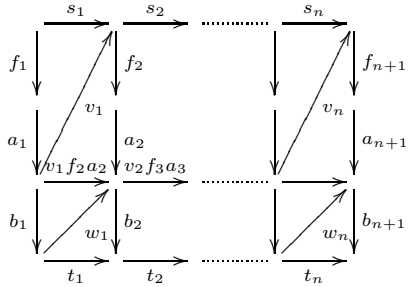
\xymatrix{
\ar[r]^{s_1}\ar[d]_{f_1}
&\ar[r]^{s_2}\ar[d]^{f_2}&\ar@{.}[r]&\ar[r]^{s_n}\ar[d]&\ar[d]^{f_{n+1}}\\
\ar[d]_{a_1}
&\ar[d]^{a_2}&&\ar[d]&\ar[d]^{a_{n+1}}\\
\ar[r]^{v_1f_2a_2}\ar[d]_{b_1}\ar[uur]_{v_1}&\ar[r]^{v_2f_3a_3}\ar[d]^{b_2}&\ar@{.}[r]&\ar[r]\ar[d]
\ar[uur]_{v_n}&\ar[d]^{b_{n+1}}\\
\ar[r]_{t_1}\ar[ur]_{w_1}&\ar[r]_{t_2}&\ar@{.}[r]&\ar[r]_{t_n}\ar[ur]_{w_n}&\\
}
$$
We prove that the there exists a morphism $h\in P_n$ such that
$h_i=f_ia_i$  which is equivalent to proving that $v_if_{i+1}a_{i+1}$ is
in $P$. We claim that $v_if_{i+1}a_{i+1}$ is the left lcm
of $f_{i+1}a_{i+1}$ and $w_i$ for $i=1,\ldots,n$: indeed it is a
common left multiple and if the left lcm was smaller then $b_i$ and $v_i$ would have a non trivial
common left divisor $y_i$ which would give a common left divisor $f_ia_iy_i$ of $s_i$ and $f_ig_i$ 
greater than their gcd $f_ia_i$. We conclude as the lcm of two morphisms in $P$ is in $P$.

So we have a morphism $h\in P_n$ which divides $fg$. There cannot be a
greater simple divisor $k$ of $fg$ as $k_i$ has to divide $s_i$
and $f_ig_i$ for $i=1\ldots n$ and $k_{n+1}$ has to be a simple divisor of
$f_{n+1}g_{n+1}$.
\end{proof}
Note that in the above proof we have used a left lcm. It is the only place where we use the fact that
$P$ is right locally Garside.
\begin{proposition}\label{alpha(f)_i}
Let $f$ be a morphism in $C(P_n)$ with source $(s_1,\ldots,s_n)$; then $\alpha(f)$ is the
simple morphism with same source such that $\alpha(f)_i= \gcd(f_i,s_i)$ for $i=1,\ldots,n$ and
$\alpha(f)_{n+1}=\alpha(f_{n+1})$.
\end{proposition}
\begin{proof}
We  write  $f=f^1\ldots  f^k$  with   $f^i\in  P_n$.  The  proof  is  by
induction  on  $k$.  The  above  lemma  proves  the  result  for  $k=2$.
We  have  $\alpha(f)=\alpha(f^1\alpha(f^2\ldots  f^k))$.  By the
induction  hypothesis applied to  $f^2\ldots f^k$, the following diagram
represents $f^1\alpha(f^2\ldots  f^k)$:
$$\xymatrix{
\ar[r]^{s_1}\ar[d]_{f_1^1}
&\ar[r]^{s_2}\ar[d]^{f_2^1}&\ar@{.}[r]&\ar[r]^{s_n}\ar[d]&\ar[d]^{f_{n+1}^1}\\
\ar[r]^{t_1}\ar[d]_{g_1}\ar[ur]&\ar[r]^{t_2}\ar[d]^{g_2}&\ar@{.}[r]&\ar[r]^{t_n}\ar[d]\ar[ur]
&\ar[d]^{g_{n+1}}\\
\ar[r]\ar[ur]&\ar[r]&\ar@{.}[r]&\ar[r]\ar[ur]&\\
}
$$
where $g_i=\gcd(t_i,(f^2\ldots f^k)_i)$ for $i\leq n$ and $g_{n+1}=\alpha((f^2\ldots f^k)_{n+1})$.
We apply now the previous lemma to  the  two term  product $f^1\alpha(f^2\ldots f^k)$.
We will be done if $\gcd(f^1_ig_i,s_i)=\gcd(f_i,s_i)$ \ie, $\gcd(f^1_i\gcd(t_i,(f^2\ldots f^k)_i),s_i)
=\gcd(f_i,s_i)$ for $i\leq n$ and $\alpha(f^1_{n+1}\alpha((f^2\ldots
f^k)_{n+1})=\alpha(f_{n+1})$. The latter is true by the properties of $\alpha$. The former
is true as the right hand side is a multiple of $f^1_i$   
so has to be the product of $f^1_i$ by a common divisor of $(f^2\ldots f^k)_i $ and $t_i$.
\end{proof}
\begin{corollary}\label{f_i divise g_i} Let $f$ and $g$
be two  morphisms in  $C(P_n)$ with same  source; then  $f\preccurlyeq g$
if  and  only  if  $f_i\preccurlyeq g_i$  for  all  $i$ 
\end{corollary}
\begin{proof} Assume that  $f_i\preccurlyeq g_i$  for  all  $i$. 
If $f$ and $g$ are in $P_n$ then we are done by lemma \ref{divisible}. In general we prove the result
by induction on the length of the normal form of $f$. We first show that
$\alpha(f)\preccurlyeq\alpha(g)$: let $(s_1,\ldots,s_n)$ be the common source of $f$ and $g$; we have
$\alpha(f)_i=\gcd(f_i,s_i)\preccurlyeq\gcd(g_i,s_i)=\alpha(g)_i$ for $i\leq n$ and 
$\alpha(f)_{n+1}=\alpha(f_{n+1})\preccurlyeq\alpha(g_{n+1})=\alpha(g)_{n+1}$, whence the result 
as $\alpha(f)$ and $\alpha(g)$ are two elements of $P_n$.
After simplifying by $\alpha(f)$ we can apply the induction hypothesis which gives
that $\omega(f)\preccurlyeq \alpha(f)\inv g$, whence $f\preccurlyeq g$.
The converse is clear.
\end{proof}
\begin{corollary}\label{f caracterise par f_i} A morphism $f\in C(P_n)$ is determined uniquely by its
source and the morphisms  $f_i\in C(P)$ \end{corollary}
\begin{proof} If
two morphisms  $f$ and $g$ have  same source and $f_i=g_i$  for all $i$,
then by  the previous  corollary they  divide each  other, so  are equal.
\end{proof}
Note that corollaries \ref{f_i divise g_i} and \ref{f caracterise par f_i}
are true for any subcategory of $C(P_n)$. Proposition \ref{alpha(f)_i}
is true in $C(P_n(\Id))$. In $C(P_n(F))$ it has to be modified as follows:
\begin{corollary}
Let $P$ be as above and $F$ be as in \ref{F};
if $f$ is a morphism in $C(P_n(F))$ with source $(s_1,\ldots,s_n)$ then $\alpha(f)$ is the
simple morphism with same source such that $\alpha(f)_i= \gcd(f_i,s_i)$ for $i=1,\ldots,n$.
\end{corollary}
\begin{proof}
It is clear that the morphism defined by these conditions is the greatest divisor of the
$\alpha(f)$ computed in $P_n$ which is in $P_n(F)$.
\end{proof}
\subsection*{More on $C(P_n(\Id))$}

We now return to the case of an only left Garside $P$ and look at $C(P_n(\Id))$.
\begin{lemma}  \label{unicity} In $P_n(\Id)$,  there is at  most one morphism
between two objects.
\end{lemma}
\begin{proof}
Indeed, if $(a_1,\ldots,a_n)\xrightarrow f(b_1,\ldots,b_n)$ in $P_n(\Id)$
then we have $f_1=1$ thus $f$ is determined by increasing induction on $i$,
using the cancellation property, by the equations $a_i=f_if'_i$ and
$b_i=f'_if_{i+1}$.
\end{proof}
Let $\ba=(a_1,\ldots,a_n)$ be an object of $P_n$ and let $(b_1,\ldots,b_n)$ be
the  normal form  of $a_1a_2\ldots  a_n$ in  $C(P)$, completed  if needed by
ones;      \ie.,     the      sequence     $b_i$      is     defined     by
$b_i=\alpha(\omega^{i-1}(a_1\ldots  a_n))$ for all $i$  (we then have $b_i=1
\Rightarrow  b_k=1 \forall  k\geq i$).  In this  situation we set $\fn(\ba)=
(b_1,\ldots,b_n)$.  Using \eg.\  inductively \ref{forme  normale de  xy} we can
always  construct  at  least  one  morphism  in  $C(P_n(\Id))$  from  $\ba$ to
$\fn(\ba)$.  In particular two objects are  in the same connected component if
and only if the product of their terms is the same.
\begin{lemma}\label{unique fleche vers fn}
In $C(P_n(\Id))$, there is a unique morphism $\ba\to\fn(\ba)$.
\end{lemma}
\begin{proof}
We first show that any morphism from $\ba$ to $\fn(\ba)$ has $\Delta_\ba$ as a
left  factor.  Let  $f=f_1f_2\ldots  f_m$  be  such  a morphism, where $f_i\in
P_n(\Id)$.  By definition of $\Delta_\ba$ we have $f_1\preccurlyeq\Delta_\ba$.
Using  left cancellability we define $g_1$ by $\Delta_\ba=f_1g_1$. Since $f_2$
and  $g_1$ both divide $\Delta_{\ba_1}$, where  $\ba_1$ is the target of $f_1$
they  have a right lcm of the form $f_2g_2=g_1h_2$. By induction on $i$ we can
extend   this  process  to  get  morphisms  $g_i,h_i\in  P_n(\Id)$  such  that
$f_ig_i=g_{i-1}h_i$.  As $\fn(\ba)$ is a final object in $C(P_n(\Id))$ we have
$g_{m+1}=\Id$ whence $f_1f_2\ldots f_m=\Delta_\ba h_2\ldots h_m$.

By induction, and using Noetherianity of $C(P_n(\Id))$, we can express any map
from  $\ba$  to  $\fn(\ba)$  as  a  (necessarily  unique)  finite  product  of
$\Delta$'s, whence the lemma.
\end{proof}
\begin{lemma} In $C(P_n(\Id))$,  there is at  most one morphism
between two objects.
\end{lemma}
\begin{proof}
Suppose there exists two distincts morphisms $f,g$ from $\ba$
to  $\bb$. The fact that one morphism exists implies that $\fn(\ba)=\fn(\bb)$.
By   composing  $f$  and  $g$  with  the  canonical  morphism  from  $\bb$  to
$\fn(\bb)=\fn(\ba)$,  we get two morphisms from $\ba$ to $\fn(\ba)$, which are
distinct  by the left cancellation property of $C(P_n(\Id))$. This contradicts
\ref{unique fleche vers fn}.
\end{proof}
\subsection*{The category $C(P_\bullet(\Id))$}
We  will  now  consider  a  category  whose  objects  can be identified to all
possible decompositions of a morphism of $C(P)$ into elements of $P$. We first
define a germ $P_\bullet(\Id)$ whose set of objects is the union of the set of
objects  of all $P_n$ for $n\geq 1$;  this germ is thus graded. For morphisms,
we  start by taking all the morphisms  of $\bigcup_n P_n(\Id)$ as morphisms of
degree   $0$.  We  will  also  add  some  morphisms  of  positive  degree.  If
$\ba=(a_1,\ldots,a_m)$ is an object of degree $m$ we denote by $\ba^{[k]}$ the
object  $(a_1,\ldots,a_m,\underbrace{1,\ldots,1}_k)$ of degree  $m+k$. Then we
add a morphism $i_{\ba,k}$ from $\ba$ to $\ba^{[k]}$ which we declare to be of
degree  $k$. We add to the germ the  products of a morphism $i_{\ba,k}$ with a
morphism  of  degree  $0$.  Finally  we  add  the  relations  (\ie. define the
following   products)  $i_{\ba,k}i_{\ba^{[k]},l}=i_{\ba,k+l}$  and,  for  each
morphism  $\ba\xrightarrow f\bb$ of degree $0$  between objects of degree $m$,
the  relations $f  i_{\bb,k}=i_{\ba,k}f^{[k]}$ where  $f^{[k]}$ is  defined by
$f^{[k]}_i=f_i$ for $i\le m$ and $f_{m+1}=\ldots=f_{m+k+1}=1$.

It  follows from these relations that any product of $i_{\ba,k}$'s and of {\em
one}  morphism of degree  $0$ is in  $P_\bullet(\Id)$, and that  a morphism in
$P_\bullet(\Id)$ is unique given its source and target (using \ref{unicity}).

The category $\cpb$ generated by $P_\bullet(\Id)$ inherits a grading.
\begin{proposition}
The category $\cpb$ is left locally Garside.
\end{proposition}
\begin{proof}
The  axioms for a germ are clear.  The locally Garside germ axiom \ref{G1}(G1)
is  also clear (in  a bounded increasing  sequence the degree becomes constant
and  we are reduced to the case of $P_n(\Id)$ and \ref{G4}(G4) is clear, using
the unicity of morphisms between two objects.

We prove now \ref{G3}(G3). Consider two maps $i_{\ba,k}f$ and $i_{\ba,l}g$. We
may  assume  that  $k\le  l$.  Then  $i_{\ba,l}f^{[l-k]}$  is  a  multiple  of
$i_{\ba,k}f$  and $i_{\ba,l}$ times the lcm of  $f^{[l-k]}$ and $g$ is the lcm
of  $i_{\ba,k}f$  and  $i_{\ba,l}g$.  Indeed  any  multiple of $i_{\ba,k}f$ of
degree  $l$ is  of the  form $i_{\ba,l}h$  where by  cancellation we must have
$f\preccurlyeq  i_{\ba^{[k]},l-k}  h$;  since  any  morphism  of  degree $l-k$
extending  $f$ must start by $f i_{\bb,l-k}=i_{\ba^{[k]},l-k}f^{[l-k]}$ (where
$\bb$ is the target of $f$) we have $f^{[l-k]}\preccurlyeq h$.

The proof of \ref{G3} (G3) is similar.
\end{proof}
\begin{remark}\label{sous-categorie pleine}
Note that $C(P_n(\Id))$ is the full subcategory of $\cpb$ obtained by
restricting the objects to those of $P^n$.
\end{remark}

We can extend \ref{unique fleche vers fn} to $\cpb$:
\begin{lemma}\label{unique fleche vers fn[k]}
For any $k$ there is a unique morphism $\ba\to\fn(\ba)^{[k]}$.
\end{lemma}
\begin{proof}
Using  the  relations  in  $C(P_\bullet(\Id))$,  any  morphism  from  $\ba$ to
$\fn(\ba)^{[k]}$  is  of  the  form  $i_{\ba,k}f$  where $f$ is a morphism of
degree  $0$ from  $\ba^{[k]}$ to  $\fn(\ba)^{[k]}$. Since  $\fn(\ba)^{[k]}$ is
clearly the normal form of $\ba^{[k]}$ we get the result by \ref{unique fleche
vers fn}, using remark \ref{sous-categorie pleine}.
\end{proof}
\begin{corollary}\label{composantes de C(P_bullet))}
Two objects of $\cpb$ are  in the same connected component if
and only if the product of their terms is the same.
\end{corollary}
\begin{proposition}
Let  $\CO$ be a  functor from $C(P_\bullet(\Id))$  to a groupoid.  Let us call
{\em  elementary isomorphism} a map of the form $\CO(f)\CO(i_{\ba',1})\inv$ or
$\CO(i_{\ba',1})\CO(f)\inv$    where    $\ba\xrightarrow   f\ba^{\prime[1]}\in
P_\bullet$ is of the form $\ba=(a_1,\ldots,a_i,a_{i+1},\ldots,a_n)\xrightarrow
f \ba^{\prime [1]}=(a_1,\ldots,a_{i-1},a_ia_{i+1},\ldots,a_n,1)$. Then all the
compositions  of elementary isomorphisms  between two objects  in the image of
$\CO$ are equal.
\end{proposition}
\begin{proof}
Given  $\ba=(a_1,\ldots,a_n)$  an  object  of  $P_n$,  and given $k\ge n$, let
$g_{\ba,k}$  be the  image by  $\CO$ of  the unique map in $C(P_\bullet(\Id))$
between  $\ba$  and  $NF(\ba)^{[k-\deg(\ba)]}$  (\cf.  \ref{unique fleche vers
fn[k]}). Then for any $k\ge n$ we have $g_{\ba,k}=\CO(f)g_{\ba^{\prime[1]},k}$
and   $\CO(i_{\ba',1})g_{\ba^{\prime[1]},k}=g_{\ba',k}$  thus  the  elementary
morphism  $\CO(f)\CO(i_{\ba',1})\inv$  between  $\CO(\ba)$  and $\CO(\ba')$ is
equal to $g_{\ba,k}g_{\ba^{\prime[1]},k}\inv$. It follows that, for $k$ larger
than  the degree of all  the objects involved, we  find by composing the above
formula  along  a  path  of  elementary  isomorphisms,  that  a composition of
elementary    isomorphisms   between    $\ba$   and    $\bb$   is   equal   to
$g_{\ba,k}g_{\bb,k}\inv$. Thus all such compositions are equal.
\end{proof}
\subsection*{An application to Deligne-Lusztig varieties}
We give an application of the last proposition to the existence of generalized
Schubert  cells associated to the elements of the braid monoid. Let $\bG$ be a
reductive group over an algebraically closed field. Let $W$ be the Weyl group,
identified to the set of orbits of $\bG$ on $\CB\times\CB$, where $\CB$ is the
variety  of Borel subgroups. Let $B^+(W)$ the corresponding Artin-Tits monoid,
and  let  $\bW$  be  the  germ  of  simple  elements of $B^+(W)$ (naturally in
bijection  with $W$), so $B^+(W)=C(\bW)$.  To an object $(\bw_1,\ldots,\bw_n)$
of $C(\bW_\bullet(\Id))$ we attach the variety
$$\CO(\bw_1,\ldots,\bw_n)=\{(\bB_1,\ldots,\bB_{n+1})\in\CB^{n+1}\mid
(\bB_i,\bB_{i+1}\})\in  \CO(w_i)\},$$ where $w_i$  is the image  of $\bw_i$ in
$W$  and $\CO(w_i)$ is  the orbit of  $\bG$ in $\CB\times\CB$ corresponding to
$w_i$. To a morphism $\bw\xrightarrow f\bv$ of $C(\bW_\bullet(\Id))$, given by
$\bw_i=f_if'_i$,   and  $\bv_i=f'_if_{i+1}$   we  associate   the  isomorphism
$\CO(f):\CO(\bw)\to\CO(\bv)$  which sends $\bB_k$ to  the unique Borel subgroup
$\bB'_k$     such    that    $(\bB_k,\bB'_k,\bB_{k+1})\in\CO(f_k,f'_k)$    and
$\bB'_{n+1}=\bB_{m+1}$. To the morphism $i_{\bw,k}$ we associate the
isomorphism which maps $(\bB_1,\ldots,\bB_{n+1})$ to
$(\bB_1,\ldots,\bB_n,\underbrace{\bB_{n+1},\ldots,\bB_{n+1}}_{k+1})$.

\begin{proposition}
$\CO$ is a functor from $C(\bW_\bullet(\Id))$ to the category of 
quasi-projective varieties with isomorphisms.
\end{proposition}
\begin{proof}
We need to check that if $\bw\xrightarrow f\bv$ and $\bv\xrightarrow g\bu$ 
are such that $f,g,fg\in\bW_\bullet(\Id)$, then $\CO(f)\CO(g)=\CO(fg)$.
This results from the fact that if $(\bB_k,\bB'_k,\bB_{k+1})\in\CO(f_k,f'_k)$,
$(\bB_{k+1},\bB'_{k+1},\bB_{k+2})\in\CO(f_{k+1},f'_{k+1})$
and $(\bB'_k,\bB''_k,\bB'_{k+1})\in\CO(g_k,g'_k)$  then
$(\bB_k,\bB''_k,\bB_{k+1})\in\CO(f_kg_k,(fg)_k')$ since $g'_k=(fg)'_kf_{k+1}$.
\end{proof}

If  $\bb=\bw_1\ldots\bw_n$, since we can pass  from any decomposition of $\bb$
into  a product of elements of $\bW$ to another by elementary isomorphisms, it
follows  that varieties associated to decompositions of the same element $\bb$
of  $B^+(W)$ are  canonically isomorphic.  Passing to  the projective limit of
these  isomorphisms,  we  can  define  a  variety  $\CO(\bb)$ associated to an
element of $B^+(W)$. Note that we could also have applied theorem
\ref{Deligne} to this situation, as Deligne did  in \cite{Deligne}.

\end{document}